\documentclass[12pt]{amsart}
\usepackage{amscd,amsmath,amssymb,enumerate,amsfonts,mathrsfs,txfonts}
\usepackage{comment}
\usepackage{hyperref}
\usepackage{xcolor}
\usepackage{tikz,pgfplots,pgf}
\usepackage{vmargin}
\usepackage[all]{xy}

\tikzset{node distance=3cm, auto}
\usetikzlibrary{calc} 
\usetikzlibrary{trees} 

\setpapersize{A4}
\setmargins{2.5cm}      
{1.5cm}                        
{16.5cm}                      
{23.42cm}                   
{10pt}                          
{1cm}                           
{0pt}                          
{2cm}

\newtheorem{theorem}{Theorem}[section]
\newtheorem{proposition}[theorem]{Proposition}
\newtheorem{definition}[theorem]{Definition}
\newtheorem{corollary}[theorem]{Corollary}
\newtheorem{lemma}[theorem]{Lemma}
\newtheorem{remark}[theorem]{Remark}
\newtheorem{remarks}[theorem]{Remarks}

\def\A{\mathcal{A}}
\def\B{\mathcal{B}}
\def\F{\mathcal{F}}
\def\G{\mathcal{G}}
\def\H{\mathcal{H}}
\def\I{\mathcal{I}}
\def\K{\mathcal{K}}
\def\W{\mathcal{W}}
\def\L{\mathcal{L}}
\def\M{\mathcal{M}}

\def\C{\mathbb{C}}
\def\D{\mathbb{D}}
\def\N{\mathbb{N}}

\def\T{\mathbb{T}}

\def\Int{\mathrm{Int}}
\def\abco{\mathrm{abco}}

\def\lin{\mathrm{lin}}
\def\rank{\mathrm{rank}}
\def\rang{\mathrm{rang}}

\def\id{\mathrm{id}}
\def\Aut{\mathrm{Aut}}

\makeatletter
\@namedef{subjclassname@2020}{\textup{2020} Mathematics Subject Classification}
\makeatother

\begin{document}

\title[Compact Bloch mappings on the complex unit disc]{Compact Bloch mappings on the complex unit disc}

\author[A. Jim{\'e}nez-Vargas]{A. Jim{\'e}nez-Vargas}
\address[A. Jim{\'e}nez-Vargas]{Departamento de Matem{\'a}ticas, Universidad de Almer{\'i}a, 04120, Almer{\'i}a, Spain}
\email{ajimenez@ual.es}

\author[D. Ruiz-Casternado]{D. Ruiz-Casternado}
\address[D. Ruiz-Casternado]{Departamento de Matem{\'a}ticas, Universidad de Almer{\'i}a, 04120, Almer{\'i}a, Spain}
\email{drc446@gmail.com}

\date{\today}

\subjclass[2010]{30H30, 46E15, 46E40, 47B38}
\keywords{Vector-valued holomorphic mapping, Bloch function, compact Bloch mapping, Bloch molecule, Bloch-free Banach space.}


\begin{abstract}
The known duality of the space of Bloch complex-valued functions on the open complex unit disc $\D$ is addressed under a new approach with the introduction of the concepts of Bloch molecules and Bloch-free Banach space of $\D$. We introduce the notion of compact Bloch mapping from $\D$ to a complex Banach space and establish its main properties: invariance by M\"obius transformations, linearization from the Bloch-free Banach space of $\D$, factorization, inclusion properties, Banach ideal property and transposition on the Bloch function space. We state Bloch versions of the classical theorems of Schauder, Gantmacher and Davis--Figiel--Johnson--Pe\l czy\'nski.
\end{abstract}
\maketitle


\section{Introduction}\label{section 1}

The study of the compactness of nonlinear mappings has been a subject of much interest as a first step to extend results from Linear Functional Analysis to the nonlinear setting. In this direction, 
various authors established some results of the classical theory such as Schauder, Gantmacher and Gantmacher--Nakamura Theorems in different nonlinear environments: for instance, R. M. Aron and M. Schottenloher \cite{AroSch-76}, J. Mujica \cite{Muj-91} and R. Ryan \cite{Rya-88} in both holomorphic and polynomial settings; J. Batt \cite{Bat-70} in the nonlinear bounded context; S. Yamamuro \cite{Yam-68} in the differentiable framework; and J. M. Sepulcre, M. Villegas-Vallecillos and the first author \cite{jsv} in the Lipschitz setting.

In this paper, we will address the study of this topic in the setting of Bloch mappings. Let $\D$ be the open complex unit disc, let $X$ be a complex Banach space and let $\H(\D,X)$ denote the space of all holomorphic mappings from $\D$ to $X$. 

The \textit{normalized Bloch space} $\widehat{\B}(\D,X)$ is the Banach space of all mappings $f\in\H(\D,X)$ with $f(0)=0$ such that   
$$
\rho_{\B}(f)=\sup\left\{(1-|z|^2)\|f'(z)\|\colon z\in\D\right\}<\infty ,
$$ 
endowed with the \textit{Bloch norm} $\rho_{\B}$. The \textit{normalized little Bloch space} $\widehat{\B}_0(\D,X)$ is the norm-closed linear subspace of $\widehat{\B}(\D,X)$ formed by all those mappings $f$ so that  
$$
\lim_{|z|\to 1^-}(1-|z|^2)\|f'(z)\|=0.
$$
To simplify, we will write $\widehat{\B}(\D)$ and $\widehat{\B}_0(\D)$ instead of $\widehat{\B}(\D,\C)$ and $\widehat{\B}_0(\D,\C)$, respectively. Note that the set of all polynomials is dense in $\widehat{\B}_0(\D)$. 

The properties of Bloch functions on $\D$ are well known (see \cite{And-85}) and they have been studied by some authors to investigate problems of geometric function theory (see \cite{AndCluPom-74, AraFisPee-85,ArrBla-03,MadMat-95,Nar-90,Zhu-07}). To paraphrase R. M. Timoney \cite{Tim-80}: \textit{one of the reasons that Bloch functions are interesting is that they can be characterized in many different ways and, so, they arise in many different contexts.} 

The properties of local compactness of holomorphic mappings on complex Banach spaces were studied in \cite{AroSch-76, Rya-88}. A holomorphic mapping $f\colon X\to Y$ is said to be \textit{locally (weakly) compact} if every point $x\in X$ has a neighborhood $U_x$ such that $f(U_x)$ is relatively (weakly) compact in $Y$. The study of holomorphic mappings on Banach spaces with a relatively (weakly) compact range was initiated by J. Mujica \cite{Muj-91}. 

The novelty of our approach in this paper lies in using the Bloch range of a function in place of its local or global range. By the \textit{Bloch range} of a function $f\colon\D\to X$, we mean the set $\left\{(1-|z|^2)f'(z)\colon z\in\D\right\}\subseteq X$. Clearly, a holomorphic mapping $f\colon\D\to X$ is Bloch if its Bloch range is a bounded subset of $X$ and this motivates us to introduce the following concepts: a holomorphic mapping $f\colon\D\to X$ is said to be \textit{(weakly) compact Bloch} if its Bloch range is a relatively (weakly) compact subset of $X$. Observe that if $X$ is finite-dimensional, then every Bloch mapping is indeed both compact and weakly compact. Among others properties, we emphasize that both linear spaces of compact Bloch mappings and weakly compact Bloch mappings from $\D$ to $X$ are \textit{M\"obius-invariant} and \textit{Banach composition Bloch ideals}.

One of the key elements in our work is the concept of \textit{Bloch atom of $\D$}, namely, the function $\gamma_z\colon\widehat{\B}(\D)\to\C$ with $z\in\D$ defined by $\gamma_z(f)=f'(z)$. These elementary Bloch molecules of $\D$ generate a norm-closed linear subspace of the dual space $\widehat{\B}(\D)^*$ of $\widehat{\B}(\D)$, denoted $\G(\D)$ and called \textit{Bloch-free Banach space of $\D$} because of a \textit{universal extension property} of this space. The elements of the linear hull generated by Bloch atoms of $\D$ will be referred to as \textit{Bloch molecules of $\D$}.

This paper is organized as follows. Section \ref{section 2} contains a complete study of the concepts introduced. It is well known that $\widehat{\B}(\D)$ is a dual Banach space \cite{AndCluPom-74, AraFisPee-85}. Section \ref{section 2'5} is devoted to the study of the duality of $\widehat{\B}(\D)$ under a different approach since we prove that $\G(\D)$ is the strongly unique predual of $\widehat{\B}(\D)$ (compare to \cite[Corollary 2]{Nar-90}). In Section \ref{section 3}, we see that $\G(\D)$ is not a mere predual of $\widehat{\B}(\D)$ but provides a procedure for linearizing the derivative of vector-valued Bloch mappings: for every complex Banach space $X$ and every mapping $f\in\widehat{\B}(\D,X)$, there exists a unique bounded linear operator $S_f\colon\G(\D)\to X$ such that $S_f\circ\Gamma=f'$, where $\Gamma\colon\D\to\G(\D)$ is the holomorphic mapping that maps $z$ to $\gamma_z$ (compare to \cite[Theorem 2.1]{Muj-91}). Furthermore, $\|S_f\|=\rho_{\B}(f)$. 

In Section \ref{section 4}, we show that the concept of compact (weakly compact, finite-rank, approximable) Bloch mapping from $\D$ to $X$ is in concordance with the corresponding property of the associated bounded linear operator $S_f$ from $\G(\D)$ to $X$. Furthermore, we apply the Davis--Figiel--Johnson--Pe\l czy\'nski theorem \cite{dfjp} to factor a weakly compact Bloch mapping through reflexive spaces. We complete this section dealing the problem as to when $\widehat{\B}(\D)$ has the approximation property. Finally, we introduce the notion of transpose mapping of a Bloch mapping, and present the Bloch analogues of the Schauder, Gantmacher and Gantmacher--Nakamura Theorems. 

We now fix some notation. $X$ and $Y$ will always denote complex Banach spaces. To indicate that $X$ and $Y$ are isometrically isomorphic, we will write $X\cong Y$. We denote by $\id_X$ the identity map of $X$. As usual, $B_X$, $S_X$ and $X^*$ denote the closed unit ball, the unit sphere and the topological dual space of $X$, respectively. For a set $A\subseteq X$, $\Int(A)$, $\abco(A)$ and $\overline{\abco}(A)$ stand for the interior, the absolutely convex hull and the norm-closed absolutely convex hull of $A$ in $X$, respectively. If $E$ and $F$ are locally convex Hausdorff spaces, $\L(E;F)$ or $\L(E,F)$ denote the linear space of all continuous linear operators from $E$ into $F$. Unless stated otherwise, if $E$ and $F$ are Banach spaces, we will understand that they are endowed with their canonical norm topology. $\T$ stands for the set of all uni-modular complex numbers.


\section{Bloch-free Banach space of the complex unit disc}\label{section 2}

Let $X$ be a Banach space and let $\Omega$ be an open subset of $\C$. Let us recall that a function $f\colon\Omega \to X$ is \textit{holomorphic} if 
$$
f'(a)=\lim_{z\to a}\frac{f(z)-f(a)}{z-a}\in X
$$
exists for all $a\in\Omega$. Let $\H(\Omega,X)$ denote the space of holomorphic mappings from $\Omega$ to $X$.

Let us recall the concept of peaking function in the setting of Bloch function spaces.

\begin{definition}
A function $f\in\widehat{\B}(\D)$ with $\rho_\B(f)\leq 1$ is said to \textit{peak at a point $z\in\D$} if $(1-|z|^2)f'(z)=1$, and for each open disc $D(z,\varepsilon)\subseteq\D$ with $\varepsilon>0$, there exists $\delta>0$ such that $(1-|w|^2)|f'(w)|\leq 1-\delta$ for all $w\in\D\backslash D(z,\varepsilon)$.
\end{definition} 

Colloquially, if $f$ peaks at $z\in\D$, we have that $(1-|w|^2)|f'(w)|$ is uniformly less than 1 when $w$ is away from $z$.

\begin{proposition}\label{peaking}
For every $z\in\D$, the function $f_z\colon\D\to\C$ defined by 
$$
f_z(w)=\frac{(1-|z|^2)w}{1-\overline{z}w}\qquad (w\in\D),
$$ 
belongs to the normalized little Bloch space $\widehat{\B}_0(\D)$ with $\rho_{\B}(f_z)=1$, and peaks at $z$.  
\end{proposition}

\begin{proof}
Note that $f_z\in\H(\D)$ with $f_z(0)=0$ and 
$$
f'_z(w)=\frac{1-|z|^2}{(1-\overline{z}w)^2}\qquad (w\in\D).
$$
Since 
$$
(1-|w|^2)|f'_z(w)|=\frac{(1-|w|^2)(1-|z|^2)}{|1-\overline{z}w|^2}=\frac{|1-\overline{z}w|^2-|z-w|^2}{|1-\overline{z}w|^2}=1-\left|\frac{z-w}{1-\overline{z}w}\right|^2\leq 1
$$
for all $w\in\D$, it follows that $f_z\in\widehat{\B}(\D)$ with $\rho_{\B}(f_z)\leq 1$. Further, $(1-|z|^2)f'_z(z)=1=\rho_\B(f_z)$. From the inequality
$$
(1-|w|^2)|f'_z(w)|=(1-|w|^2)\frac{1-|z|^2}{|1-\overline{z}w|^2}\leq (1-|w|^2)\frac{1}{|1-|z||w||^2}\qquad (w\in\D),
$$
we infer that 
$$
\lim_{|w|\to 1^-}(1-|w|^2)\left|f_z'(w)\right|=0,
$$
and thus $f_z\in\widehat{\B}_0(\D)$. Finally, let $0<\varepsilon<1$. If $|w|<1$ and $|w-z|\geq\varepsilon$, we have 
$$
(1-|w|^2)|f'_z(w)|=1-\left|\frac{z-w}{1-\overline{z}w}\right|^2\leq 1-\frac{\varepsilon^2}{4},
$$
and so $f_z$ peaks at $z$.
\end{proof}

The first step to construct a Bloch molecule on the normalized Bloch space $\widehat{\B}(\D)$ is the following concept.  

\begin{definition}
For each $z\in\D$, a \textit{Bloch atom} of $\D$ is the function $\gamma_z\colon\widehat{\B}(\D)\to\C$ given by 
$$
\gamma_z(f)=f'(z)\qquad (f\in\widehat{\B}(\D)).
$$
\end{definition}

We now show that $\gamma_z$ is an element of the topological dual space of $\widehat{\B}(\D)$.

\begin{proposition}\label{p1}
Let $z\in\D$. Then $\gamma_z$ is in $\widehat{\B}(\D)^*$ with $\|\gamma_z\|=1/(1-|z|^2)$.
\end{proposition}

\begin{proof}
Clearly, $\gamma_z$ is linear and 
$$
\left|\gamma_z(f)\right|=\left|f'(z)\right|\leq\frac{\rho_\B(f)}{1-|z|^2}
$$
for all $f\in\widehat{\B}(\D)$. Hence $\gamma_z$ is continuous with $\|\gamma_z\|\leq 1/(1-|z|^2)$. To see that $\|\gamma_z\|=1/(1-|z|^2)$, take the function $f_z\in\widehat{\B}(\D)$. Since $\rho_{\B}(f_z)=1$, we have 
$$
\|\gamma_z\|\geq \left|\gamma_z(f_z)\right|=\left|f'_z(z)\right|=\frac{1}{1-|z|^2}.
$$
\end{proof}

The preceding proposition justifies the following concept.

\begin{definition}
A \textit{normalized Bloch atom} of $\D$ is the function $\widehat{\gamma}_z\colon\widehat{\B}(\D)\to\C$ given by
$$
\widehat{\gamma}_z=(1-|z|^2)\gamma_z.
$$
We denote by $\M_\B(\D)$ the \textit{set of all normalized Bloch atoms of $\D$}, that is, 
$$
\M_\B(\D)=\left\{\widehat{\gamma}_z\colon z\in\D\right\}.
$$
\end{definition}

We now introduce the concept of Bloch-free Banach space of $\D$. This terminology will be justified by a \textit{universal extension property} of this space (see Theorem \ref{teo-3}).

\begin{definition}
The \textit{Bloch-free Banach space over} $\D$, denoted $\G(\D)$, is the norm-closed linear hull of $\left\{\gamma_z\colon z\in\D\right\}$ in $\widehat{\B}(\D)^*$, that is, 
$$
\G(\D)=\overline{\lin}\left\{\gamma_z\colon z\in\D\right\}\subseteq\widehat{\B}(\D)^*.
$$
The elements of $\lin\{\gamma_z\colon z\in\D\}\subseteq\widehat{\B}(\D)^*$ are referred to as \textit{Bloch molecules} of $\D$. 
\end{definition}

In consistency with the notation of Bloch atoms, we consider the mapping $\Gamma\colon\D\to\widehat{\B}(\D)^*$ defined by 
$$
\Gamma(z)=\gamma_z\qquad (z\in\D).
$$
In order to prove its holomorphy, let us recall that given a Banach space $X$, a subset $N\subseteq B_{X^*}$ is said to be \textit{norming for $X$} if  
$$
\|x\|=\sup\left\{\left|x^*(x)\right|\colon x^*\in N\right\}\qquad (x\in X).
$$
Notice that $\M_\B(\D)$ is norming for $\widehat{\B}(\D)$ since
$$
\rho_{\B}(f)=\sup\left\{(1-|z|^2)\left|f'(z)\right|\colon z\in\D\right\}=\sup\left\{\left|\widehat{\gamma}_z(f)\right|\colon z\in\D\right\}
$$
for every $f\in\widehat{\B}(\D)$. A mapping $f\colon\D\to X$ is said to be \textit{locally bounded} if 
$$
\sup\left\{\|f(z)\|\colon z\in K\right\}<\infty
$$
for every compact set $K\subseteq\D$. 

\begin{proposition}\label{prop-1}
$\Gamma$ is a holomorphic mapping from $\D$ to $\G(\D)$ with $\Gamma'(z)(f)=f''(z)$ for all $f\in\widehat{\B}(\D)$ and $z\in\D$.
\end{proposition}

\begin{proof}
Start by noting that for each $f\in\widehat{\B}(\D)$, the \textit{evaluation mapping} $J_f\colon\G(\D)\to\C$, defined by 
$$
J_f(\gamma)=\gamma(f)\qquad (\gamma\in\G(\D)),
$$
is linear and continuous with $\|J_f\|\leq \rho_\B(f)$. It is clear that $\{J_f\colon f\in B_{\widehat{\B}(\D)}\}$ is a norming set for $\G(\D)$ since $\|\gamma\|=\sup\{\left|J_f(\gamma)\right|\colon f\in B_{\widehat{\B}(\D)}\}$ for each $\gamma\in\G(\D)$. 
On the other hand, given a compact set $K\subseteq\D$, there exists $r_K\in (0,1)$ such that $|z|\leq r_K$ for all $z\in K$. 
Using Proposition \ref{p1}, we have   
$$
\|\Gamma(z)\|=\|\gamma_z\|=\frac{1}{1-|z|^2}\leq\frac{1}{1-(r_K)^2}
$$
for all $z\in K$, and thus 
$$
\sup\left\{\|\Gamma(z)\|\colon z\in K\right\}\leq\frac{1}{1-(r_K)^2}.
$$
Hence $\Gamma\colon\D\to\G(\D)$ is locally bounded. Since $J_f\circ\Gamma=f'\in\H(\D)$ for all $f\in B_{\B(\D)}$, then $\Gamma$ is holomorphic by applying \cite[Proposition A.3]{AreBatHieNeu-01}. 

Finally, given $a\in\D$ and $f\in\widehat{\B}(\D)$, we have
$$
\left(\frac{\Gamma(z)-\Gamma(a)}{z-a}\right)(f)=\frac{\Gamma(z)(f)-\Gamma(a)(f)}{z-a}=\frac{f'(z)-f'(a)}{z-a}\qquad (z\in\D,\, z\neq a),
$$
hence
$$
\lim_{z\to a}\left(\frac{\Gamma(z)-\Gamma(a)}{z-a}\right)(f)=f''(a),
$$
and thus $\Gamma'(a)(f)=f''(a)$. 
\end{proof}

We gather some properties of $\G(\D)$ in the following.

\begin{remarks}\label{remark1}
\begin{enumerate}
\item Since $\D$ is separable and the mapping $\Gamma\colon\D\to\G(\D)$ is continuous, the set $\Gamma(\D)$ is separable, hence its norm-closed linear hull, $\G(\D)$, is also separable. 
\item $\Gamma(\D)$ is a linearly independent subset of $\G(\D)$. Indeed, assume that $\sum_{k=1}^n\lambda_k\gamma_{z_k}=0$, where $n\in\N$, $\lambda_1,\ldots,\lambda_n\in\C$ and $z_1,\ldots,z_n\in\D$ with $z_j\neq z_k$ if $j,k\in\{1,\ldots,n\}$ and $j\neq k$. Then $\sum_{k=1}^n\lambda_kf'(z_k)=0$ for all $f\in\widehat{\B}(\D)$. For each $j\in\{1,\ldots,n\}$, let $Q_j\colon\D\to\C$ be the polynomial defined by 
$$
Q_j(z)=\prod_{k=1,\\ k\neq j}^n\frac{z-z_k}{z_j-z_k},
$$
and let $P_j\colon\D\to\C$ be a polynomial such that $P'_j=Q_j$ and $P_j(0)=0$. Then $P_j$ is in $\widehat{\B}(\D)$ with $P'_j(z_k)=\delta_{jk}$ for all $k\in\{1,\ldots,n\}$, where $\delta$ is the Kronecker delta. It follows that  
$$
0=\sum_{k=1}^n\lambda_k P'_j(z_k)=\sum_{k=1}^n\lambda_k \delta_{jk}=\lambda_j
$$
for each $j\in\{1,\ldots,n\}$, as we wanted.
\end{enumerate}
\end{remarks}

We also will use the following easy fact through the paper.

\begin{lemma}\label{main lemma}
Let $X$ be a complex Banach space and let $F\in\mathcal{H}(\mathbb{D},X)$. Then there exists a mapping $f\in\mathcal{H}(\mathbb{D},X)$ with $f(0)=0$ such that $f'=F$.
\end{lemma} 

\begin{proof}
Define $f\colon\mathbb{D}\to X$ by 
$$
f(z)=\int_0^z F(s)\ ds\qquad (z\in\mathbb{D}),
$$
where the integral is the Bochner integral of the function $F\colon\mathbb{D}\to X$ along some path in $\mathbb{D}$ from $0$ to $z$. First, notice that $f$ does not depend on the path chosen since if $\alpha_z,\beta_z$ are paths in $\mathbb{D}$ from $0$ to $z$ and $x^*\in X^*$, by classical Cauchy's Theorem for star-shaped domains we have 
$$
\int_{\alpha_z}(x^*\circ F)(s)\ ds=\int_{\beta_z}(x^*\circ F)(s)\ ds,
$$
that is,
$$
x^*\left(\int_{\alpha_z}F(s)\ ds\right)=x^*\left(\int_{\beta_z}F(s)\ ds\right),
$$ 
by the properties of the Bochner integral, and since $X^*$ separates the points of $X$, it follows that 
$$
\int_{\alpha_z}F(s)\ ds=\int_{\beta_z}F(s)\ ds.
$$
Clearly, $f(0)=0$. Now, given $z\in\mathbb{D}$ and $x^*\in X^*$, we have  
$$
(x^*\circ f)(z)=x^*\left(\int_0^z F(s)\ ds\right)=\int_0^z (x^*\circ F)(s)\ ds,
$$
and again Cauchy's Theorem for star-shaped domains gives
$$
x^*(f'(z))=(x^*\circ f)'(z)=(x^*\circ F)(z)=x^*(F(z)),
$$
and this yields that $f'(z)=F(z)$. 
\end{proof}


\section{Duality of complex-valued Bloch function spaces}\label{section 2'5}

The following theorem shows that $\G(\D)$ is a predual of $\widehat{\B}(\D)$.

\begin{theorem}\label{teo-2}
The space $\widehat{\B}(\D)$ is isometrically isomorphic to $\G(\D)^*$, via the mapping $\Lambda\colon\widehat{\B}(\D)\to\G(\D)^*$ given by
$$
\Lambda(f)(\gamma)=\sum_{k=1}^n\lambda_kf'(z_k)\qquad (f\in\widehat{\B}(\D),\; \gamma=\sum_{k=1}^n\lambda_k\gamma_{z_k}\in\lin(\Gamma(\D))).
$$
Its inverse is the mapping $\Lambda^{-1}\colon\G(\D)^*\to\widehat{\B}(\D)$ defined by 
$$
\Lambda^{-1}(\phi)(z)=\int_{0}^z\phi(\gamma_w)\ dw \qquad \left(\phi\in\G(\D)^*,\; z\in\D\right).
$$
\end{theorem}

\begin{proof}
Let $f\in\widehat{\B}(\D)$. Define the mapping $\Lambda_0(f)\colon\lin(\Gamma(\D))\to\C$ by 
$$
\Lambda_0(f)(\gamma)=\sum_{k=1}^n\lambda_k f'(z_k)\qquad (\gamma=\sum_{k=1}^n\lambda_k\gamma_{z_k}).
$$
Clearly, $\Lambda_0(f)$ is linear. Given $\gamma=\sum_{k=1}^n\lambda_k\gamma_{z_k}\in\lin(\Gamma(\D))$, we have  
$$
\left|\Lambda_0(f)(\gamma)\right|=\left|\gamma(f)\right|\leq \rho_{\B}(f)\|\gamma\|.
$$
Hence $\Lambda_0(f)$ is continuous with $\|\Lambda_0(f)\|\leq \rho_{\B}(f)$. In fact, $\|\Lambda_0(f)\|=\rho_{\B}(f)$ since 
$$
\|\Lambda_0(f)\| \geq \sup_{z\in\mathbb{D}} |\Lambda_0(f)(\widehat{\gamma}_z)| = \sup_{z\in\mathbb{D}} (1-|z|^2)|f'(z)| =\rho_{\mathcal{B}}(f).
$$
By the denseness of $\lin(\Gamma(\D))$ in $\G(\D)$, there exists a unique continuous function $\Lambda(f)$ from $\G(\D)$ to $\C$ that extends $\Lambda_0(f)$. Further, $\Lambda(f)$ is linear and $\|\Lambda(f)\|=\|\Lambda_0(f)\|$. 
Let $\Lambda\colon\widehat{\B}(\D)\to\G(\D)^*$ be the mapping so defined. Clearly, $\Lambda$ is linear. In order to see that it is a surjective isometry, let $\phi$ be an element of $\G(\D)^*$. Define $F\colon\D\to\C$ by
$$
F(z)=\phi(\gamma_z)\qquad\left(z\in\D\right).
$$ 
Notice that $F=\phi\circ\Gamma\in\H(\D)$.  
Now, by Lemma \ref{main lemma}, there exists a holomorphic function $f\colon\D\to\C$ with $f(0)=0$ such that $f'(z)=\phi(\gamma_z)$ for all $z\in\D$. Since 
$$
(1-|z|^2)|f'(z)|=(1-|z|^2)\left|\phi(\gamma_z)\right|\leq (1-|z|^2)\|\phi\|\|\gamma_z\|=\|\phi\|
$$
for all $z\in\D$, we deduce that $f\in\widehat{\B}(\D)$ with $\rho_{\B}(f)\leq\|\phi\|$. For any $\gamma=\sum_{k=1}^n\lambda_k\gamma_{z_k}\in\lin(\Gamma(\D))$, we get 
$$
\Lambda(f)(\gamma)
=\Lambda_0(f)(\gamma)
=\sum_{k=1}^n\lambda_k f'(z_k)
=\sum_{k=1}^n\lambda_k\phi(\gamma_{z_k})
=\phi\left(\sum_{k=1}^n\lambda_k\gamma_{z_k}\right)
=\phi(\gamma).
$$
Hence $\Lambda(f)=\phi$ on the dense subspace $\lin(\Gamma(\D))$ of $\G(\D)$ and, consequently, $\Lambda(f)=\phi$. This completes the proof of the theorem. 
\end{proof}

We can now use the classic notation of the dual pairing between $\widehat{\B}(\D)$ and $\G(\D)$:
$$
\left\langle \gamma,f\right\rangle=\gamma(f)\qquad (\gamma\in\G(\D),\; f\in \widehat{\B}(\D)).
$$


\begin{corollary}\label{remark1-1}
The closed unit ball of $\G(\D)$ coincides with the norm-closed absolutely convex hull of $\M_\B(\D)$ in $\widehat{\B}(\D)^*$.
\end{corollary}

\begin{proof}
Since $\M_\B(\D)\subseteq S_{\G(\D)}$, it is immediate that $\overline{\abco}(\M_\B(\D))\subseteq B_{\G(\D)}$. 
Conversely, suppose that there is a $\upsilon\in B_{\G(\D)}\setminus\overline{\abco}(\M_\B(\D))$. By the Bipolar Theorem, 
there exists a $\phi\in\G(\D)^*$ such that $\left|\phi(\gamma)\right|\leq 1<\left|\phi(\upsilon)\right|$ for all $\gamma\in\abco(\M_\B(\D))$. By Theorem \ref{teo-2}, $\phi=\Lambda(f)$ for some $f\in\widehat{\B}(\D)$. It follows that
$$
(1-|z|^2)\left|f'(z)\right|=\left|\Lambda(f)(\widehat{\gamma}_z)\right|\leq 1<\left|\Lambda(f)(\upsilon)\right|\leq \|\Lambda(f)\|\|\upsilon\|\leq \|\Lambda(f)\|=\rho_\B(f)
$$
for all $z\in\D$, and we arrive at a contradiction. Hence $B_{\G(\D)}\subseteq\overline{\abco}(\M_\B(\D))$, as required.
\end{proof}

We have just seen that $B_{\G(\D)}=\overline{\abco}(\M_\B(\D))$. In fact, a stronger approximation result is valid: every element of $\G(\D)$ can be represented as a series of normalized Bloch atoms, and one may choose the sum of the coefficients to be arbitrarily close to the norm.  

\begin{theorem}
Let $\gamma\in\G(\D)$. Then, for every $\varepsilon>0$, there exist sequences $(\lambda_n)_{n\in\N}\in\ell_1$ with $\sum_{n=1}^\infty\left|\lambda_n\right|<\|\gamma\|+\varepsilon$ and $(z_n)_{n\in\N}\in\D^\N$ such that  
$$
\gamma=\sum_{n=1}^\infty \lambda_n\widehat{\gamma}_{z_n}.
$$
Moreover,
$$
\|\gamma\|=\inf\left\{\sum_{n=1}^\infty\left|\lambda_n\right|\colon \gamma=\sum_{n=1}^\infty \lambda_n\widehat{\gamma}_{z_n},\ (\lambda_n)_{n\in\N}\in\ell_1,\; (z_n)_{n\in\N}\in\D^\N\right\}.
$$
The dual pairing satisfies the formula
$$
\left\langle\sum_{n=1}^\infty \lambda_n\widehat{\gamma}_{z_n},f\right\rangle=\sum_{n=1}^\infty \lambda_n(1-|z_n|^2)f'(z_n)\qquad \left(f\in\widehat{\B}(\D)\right).
$$
\end{theorem}

\begin{proof} 
Let $A=\left\{z_n\colon n\in\N\right\}$ be a dense countable subset of $\D$ and consider the normed space $\ell_1(A)$ of all functions $f\colon A\to\C$ such that $\sum_{n=1}^{\infty}|f(z_n)|<\infty$ with the norm
$$
\|f\|_1=\sum_{n=1}^{\infty}|f(z_n)|.
$$
For each $n\in\N$, let $e_{z_n}\colon A\to\C$ be the function defined by $e_{z_n}(w)=\chi_{z_nw}$, where $\chi$ is the Kronecker delta. Clearly, $e_{z_n}\in\ell_1(A)$ with $\|e_{z_n}\|_1=1$. 
Given $f\in\ell_1(A)$, the series $\sum_{n\in\N}f(z_n)\widehat{\gamma}_{z_n}$ is absolutely convergent since  
\[
\sum _{k=1}^n\|f(z_k)\widehat{\gamma}_{z_k}\|=\sum_{k=1}^n\left|f(z_k)\right|\|\widehat{\gamma}_{z_k}\|=\sum _{k=1}^n\left|f(z_k)\right|\qquad (n\in \N). 
\]
Since $\G(\D)$ is a Banach space, then the series $\sum_{n\in\N}f(z_n)\widehat{\gamma}_{z_n}$ converges in $\G(\D)$ and therefore we can consider the mapping $T\colon\ell_1(A)\to\G(\D)$ defined by  
\[
T(f)=\sum_{n=1}^\infty f(z_n)\widehat{\gamma}_{z_n}\qquad \left(f=\sum_{n=1}^\infty f(z_n)e_{z_n}\in\ell_1(A)\right). 
\]
Clearly, $T\in\L(\ell_1(A),\G(\D))$ with $\|T\|\leq 1$. Now, note that $T(e_{z_n})=\widehat{\gamma}_{z_n}$ for each $n\in\N$
and therefore
$$
B_{\G(\D)}=\overline{\abco}(\M_\B(\D))\subseteq\overline{T(B_{\ell_1(A)})}.
$$
By a standard Open Mapping Theorem argument (see Lemma 2.23 and the remark which follows it in \cite{fhhmz}), we actually have 
$$
\Int(B_{\G(\D)})\subseteq T(\Int(B_{\ell_1(A)})).
$$
Let $\gamma\in\G(\D)$ and $\varepsilon>0$. Since $\gamma/(\|\gamma\|+\varepsilon)\in\Int(B_{\G(\D)})$, there exists $f\in\ell_1(A)$ with $\sum_{n=1}^{\infty}|f(z_n)|<1$ such that 
$$
\gamma/(\|\gamma\|+\varepsilon)=T(f)=\sum_{n=1}^\infty f(z_n)\widehat{\gamma}_{z_n}.
$$
Take $\lambda_n=(\|\gamma\|+\varepsilon)f(z_n)$ for all $n\in\N$, and the first assertion is proved. The second assertion follows from the first one and the fact that 
$$
\|\gamma\|\leq\sum_{n=1}^\infty\|\lambda_n\widehat{\gamma}_{z_n}\|=\sum_{n=1}^\infty\left|\lambda_n\right|
$$
for any such expression of $\gamma$ as in the statement. 

To check the dual pairing formula, consider such a representation of $\gamma$, $\sum_{n=1}^\infty \lambda_n\widehat{\gamma}_{z_n}$. Given $f\in\widehat{\B}(\D)$, we must see that the series $\sum_{n\in\N} \lambda_n(1-|z_n|^2)f'(z_n)$ converges to $\gamma(f)$. Denote, for each $n\in \mathbb{N}$, $\gamma_n=\sum_{k=1}^n \lambda_k\widehat{\gamma}_{z_k}\in\lin(\M_\B(\D))$. Since $\sum_{k=1}^n \lambda_k(1-|z_k|^2)f'(z_k)=\gamma_n(f)$, we have
$$
\left|\gamma(f)-\sum_{k=1}^n\lambda_k(1-|z_k|^2)f'(z_k)\right|=
|\gamma(f)-\gamma_n(f)|\leq \|\gamma-\gamma_n\|\rho_\B(f)\qquad (n\in\N).
$$
Since $\|\gamma-\gamma_n\|\to 0$, we get $\sum_{n=1}^\infty\lambda_n(1-|z_n|^2)f'(z_n)=\gamma(f)$.
\end{proof}

The isometric isomorphism $\Lambda\colon\widehat{\B}(\D)\to\G(\D)^*$ in Theorem \ref{teo-2} allows to consider the weak* topology on $\widehat{\B}(\D)$, that is, the topology
$$
\left\{\Lambda^{-1}(U)\colon U \text{ is open in } \left(\G(\D)^*,w^*\right)\right\}.
$$
On bounded subsets of $\widehat{\B}(\D)$, the weak* topology on $\widehat{\B}(\D)$ agrees with the topology of pointwise convergence as we see in the next result whose proof is straightforward.

\begin{corollary}\label{cor-1}
Let $(f_i)_{i\in I}$ be a net in $\widehat{\B}(\D)$ and $f\in\widehat{\B}(\D)$.
\begin{enumerate}
	\item If $(f_i)_{i\in I}\to f$ weak* in $\widehat{\B}(\D)$, then $(f_i)_{i\in I}\to f$ pointwise on $\D$.
	\item If $(f_i)_{i\in I}$ is bounded in $\widehat{\B}(\D)$ and $(f_i)_{i\in I}\to f$ pointwise on $\D$, then $(f_i)_{i\in I}\to f$ weak* in $\widehat{\B}(\D)$. $\hfill\Box$
\end{enumerate}
\end{corollary}

If $h\colon\D\to\D$ is a holomorphic function, then the Pick--Schwarz Lemma asserts that  
$$
(1-|z|^2)|h'(z)|\leq 1-|h(z)|^2\qquad (z\in\D).
$$
Let $h\colon\D\to \D$ be a holomorphic function with $h(0)=0$ and $f\in\widehat{\B}(\D)$. Clearly, $(f\circ h)(0)=0$ and 
$$
(1-|z|^2)|(f\circ h)'(z)|=(1-|z|^2)|f'(h(z))|h'(z)|\leq(1-|h(z)|^2)|f'(h(z))|\leq \rho_\B(f)
$$
for all $z\in\D$, and thus $f\circ h\in\widehat{\B}(\D)$ with $\rho_\B(f\circ h)\leq \rho_\B(f)$. So we have the following.

\begin{proposition}
For each holomorphic function $h\colon\D\to\D$ with $h(0)=0$, the composition operator $C_h\colon\widehat{\B}(\D)\to\widehat{\B}(\D)$, defined by $C_h(f)=f\circ h$ for all $f\in\widehat{\B}(\D)$, is linear and continuous with $\|C_h\|\leq 1$. Moreover, $\|C_h\|\geq \rho_\B(C_h(\id_\D))=\rho_\B(h)$. $\hfill\Box$
\end{proposition}

\begin{corollary}\label{cor-3}
Let $h\colon\D\to \D$ be a holomorphic function with $h(0)=0$. There exists a unique operator $\widehat{h}\in\L(\G(\D),\G(\D))$ such that $\widehat{h}\circ\Gamma=h'\cdot(\Gamma\circ h)$. In fact, $(\widehat{h})^*=\Lambda\circ C_h\circ\Lambda^{-1}$ and therefore $\|\widehat{h}\|=\|C_h\|$. 
\end{corollary}

\begin{proof}
First, we claim that $\Lambda\circ C_h\circ\Lambda^{-1}$ is weak*-to-weak* continuous from $\G(\D)^*$ into itself. Let $(\phi_i)_{i\in I}$ be a net in $\G(\D)^*$ and $\phi\in\G(\D)^*$. Assume that $(\phi_i)_{i\in I}\to\phi$ weak* in $\G(\D)^*$. By \cite[Corollary 2.6.10]{Meg-98}, $(\phi_i)_{i\in I}$ is norm-bounded in $\G(\D)^*$. Clearly, $(\Lambda^{-1}(\phi_i))_{i\in I}\to\Lambda^{-1}(\phi)$ weak* in $\widehat{\B}(\D)$. Hence $(\Lambda^{-1}(\phi_i))_{i\in I}\to\Lambda^{-1}(\phi)$ pointwise on $\D$ by Corollary \ref{cor-1} (1). In particular, $(\Lambda^{-1}(\phi_i)\circ h)_{i\in I}\to\Lambda^{-1}(\phi)\circ h$ pointwise on $\D$. Since $(\Lambda^{-1}(\phi_i)\circ h)_{i\in I}$ is norm-bounded because 
$$
\rho_\B(\Lambda^{-1}(\phi_i)\circ h)\leq \rho_\B(\Lambda^{-1}(\phi_i))=\|\Lambda(\Lambda^{-1}(\phi_i))\|=\|\phi_i\|
$$ 
for all $i\in I$, Corollary \ref{cor-1} (2) guarantees that $(\Lambda^{-1}(\phi_i)\circ h)_{i\in I}\to\Lambda^{-1}(\phi)\circ h$ weak* in $\widehat{\B}(\D)$, that is, $(C_h(\Lambda^{-1}(\phi_i)))_{i\in I}\to C_h(\Lambda^{-1}(\phi))$ weak* in $\widehat{\B}(\D)$. Finally, $(\Lambda(C_h(\Lambda^{-1}(\phi_i))))_{i\in I}\to \Lambda(C_h(\Lambda^{-1}(\phi)))$ weak* in $\G(\D)^*$ and this proves our claim. 

By \cite[Corollaries 3.1.11 and 3.1.5]{Meg-98}, there is a unique operator $\widehat{h}\in\L(\G(\D),\G(\D))$ such that $(\widehat{h})^*=\Lambda\circ C_h\circ\Lambda^{-1}$. It is clear that $\|\widehat{h}\|
=\|(\widehat{h})^*\|=\|C_h\|$. Given $f\in\widehat{\B}(\D)$ and $z\in\D$, we have
$$
(\widehat{h})^*(\Lambda(f))(\gamma_z)=(\Lambda(f)\circ\widehat{h})(\gamma_z)=\Lambda(f)(\widehat{h}(\gamma_z)),
$$
and
\begin{align*}
\Lambda\circ C_h\circ\Lambda^{-1}(\Lambda(f))(\gamma_z)&=\Lambda(C_h(f))(\gamma_z)=\Lambda(f\circ h)(\gamma_z)\\
																												&=(f\circ h)'(z)=f'(h(z))h'(z)\\
																												&=\Lambda(f)(\gamma_{h(z)})h'(z)=\Lambda(f)(h'(z)\gamma_{h(z)}).
\end{align*}
Hence 
$$
\Lambda(f)(\widehat{h}(\gamma_z))=\Lambda(f)(h'(z)\gamma_{h(z)})
$$
for all $f\in\widehat{\B}(\D)$ and $z\in\D$. By the surjectivity of $\Lambda\colon\widehat{\B}(\D)\to\G(\D)^*$, this means that  
$$
\phi(\widehat{h}(\gamma_z))=\phi(h'(z)\gamma_{h(z)})
$$
for all $\phi\in\G(\D)^*$ and $z\in\D$. Since $\G(\D)^*$ separates points of $\G(\D)$, this implies that $\widehat{h}(\gamma_z)=h'(z)\gamma_{h(z)}$ for all $z\in\D$, that is, $\widehat{h}\circ\Gamma=h'\cdot(\Gamma\circ h)$.   
\end{proof}

Let us recall (see \cite[B.3]{Pie-80}) that if $X$ and $Y$ are Banach spaces and $T\in\L(X,Y)$, then $T$ is called a \emph{metric injection or an isometry} if $\|T(x)\|=\|x\|$ for all $x\in X$, and $T$ is called a \emph{metric surjection} if $T$ is surjective and $\|T(x)\|=\inf\left\{\|y\|\colon T(y)=T(x)\right\}$ for all $x\in X$.

\begin{corollary}\label{cor-august}
Let $h\colon\D\to\D$ be a holomorphic function such that $h(0)=0$.
\begin{enumerate}
	\item If $g\colon\D\to\D$ is a holomorphic function with $g(0)=0$, then $\widehat{g\circ h}=\widehat{g}\circ\widehat{h}$.
	\item $\widehat{\id_\D}=\id_{\G(\D)}$.
	\item If $h$ is bijective and $h^{-1}$ is holomorphic, then $\widehat{h}$ is a topological isomorphism from $\G(\D)$ onto itself with $(\widehat{h})^{-1}=\widehat{h^{-1}}$.
	\item $\widehat{h}\colon\G(\D)\to\G(\D)$ is a metric surjection if and only if $C_h\colon\widehat{\B}(\D)\to\widehat{\B}(\D)$ is a metric injection.
	\item $\widehat{h}\colon\G(\D)\to\G(\D)$ is a metric injection if and only if $C_h\colon\widehat{\B}(\D)\to\widehat{\B}(\D)$ is a metric surjection.
\end{enumerate}
\end{corollary}

\begin{proof}
(1) By Corollary \ref{cor-3}, there are unique operators $\widehat{h},\widehat{g}\in\L(\G(\D),\G(\D))$ such that $\widehat{h}\circ\Gamma=h'\cdot(\Gamma\circ h)$ and $\widehat{g}\circ\Gamma=g'\cdot(\Gamma\circ g)$. Since 
\begin{align*}
(\widehat{g}\circ\widehat{h})\circ\Gamma&=\widehat{g}\circ(\widehat{h}\circ\Gamma)=\widehat{g}\circ(h'\cdot(\Gamma\circ h))\\
																				&=h'\cdot(\widehat{g}\circ(\Gamma\circ h))=h'\cdot((\widehat{g}\circ\Gamma)\circ h)\\
																				&=h'\cdot((g'\cdot(\Gamma\circ g))\circ h)=h'\cdot((g'\circ h)\cdot((\Gamma\circ g)\circ h))\\
																				&=(g\circ h)'\cdot(\Gamma\circ(g\circ h))
\end{align*}
it follows that $\widehat{g\circ h}=\widehat{g}\circ\widehat{h}$.

(2) Since $\id_{\G(\D)}\circ\Gamma=\Gamma=(\id_\D)'\cdot(\Gamma\circ\id_\D)$, we deduce that $\widehat{\id_\D}=\id_{\G(\D)}$. 

(3) It is easily obtained from Corollary \ref{cor-3} and parts (1) and (2).

(4)-(5) Since $(\widehat{h})^*$ is $C_h$ up to isometric isomorphisms, it follows from \cite[B.3.9]{Pie-80}.
\end{proof}

In the next result, we will make of use the following auxiliary functions. Given $f\in\widehat{\B}(\D)$ and $r\in (0,1)$, it is easy to prove that the function $f_r(z)=f(rz)$ for all $z\in\D$ belongs to $\widehat{\B}_0(\D)$ with $\rho_{\B}(f_r)\leq \rho_\B(f)$. 

Applying Corollary \ref{cor-1}, we deduce the following property of normalized little Bloch space.

\begin{corollary}\label{cor-2}
$\widehat{\B}_0(\D)$ is weak*-dense in $\widehat{\B}(\D)$.
\end{corollary}

\begin{proof}
Let $f\in\widehat{\B}(\D)$ and let $(r_n)_{n\in\N}$ be a sequence in $(0,1)$ converging to $1$. For each $n\in\mathbb{N}$, the function $f_{r_n}\in\widehat{\B}_0(\D)$ with $\rho_\B(f_{r_n})\leq \rho_\B(f)$, and $f_{r_n}(z)=f(r_nz)\to f(z)$ as $n\to +\infty$ for every $z\in\D$. By Corollary \ref{cor-1}, $(f_{r_n})_{n\in\N}$ converges weak* to $f$. This proves the corollary. 
\end{proof}

We will now prove that $\widehat{\B}_0(\D)$ is a predual of $\G(\D)$. Let $\mathcal{C}_b(\D)$ be the Banach space of all bounded continuous complex-valued functions on $\D$ with the supremum norm, and let $\mathcal{C}_0(\D)$ be its closed subspace of functions that vanish at infinity. We begin noting that the mapping $\Phi\colon\widehat{\B}(\D)\to\mathcal{C}_b(\D)$, defined by 
$$
\Phi(f)(z)=(1-|z|^2)f'(z) \qquad (f\in\widehat{\B}(\D), \; z\in\D),
$$
is an isometric linear embedding from $\widehat{\B}(\D)$ into $\mathcal{C}_b(\D)$. 

\begin{theorem}\label{messi1}
The restriction mapping $R\colon\G(\D)\to\widehat{\B}_0(\D)^*$, defined by
$$
R(\gamma)(f)=\gamma(f) \qquad \left(f\in\widehat{\B}_0(\D), \; \gamma\in\G(\D)\right),
$$
is an isometric isomorphism. 
\end{theorem}

\begin{proof}
Since $\G(\D)\subseteq\widehat{\B}(\D)^*$, it is clear that $R$ is a linear mapping from $\G(\D)$ into $\widehat{\B}_0(\D)^*$ satisfying $\|R(\gamma)\|\leq\|\gamma\|$ for all $\gamma\in\G(\D)$. We next prove that $R$ is surjective. Take $\varphi\in\widehat{\B}_0(\D)^*$. The functional $T\colon\Phi(\widehat{\B}_0(\D))\to\C$, defined by $T(\Phi(f))=\varphi(f)$ for all $f\in\widehat{\B}_0(\D)$, is linear, continuous and $\|T\|=\|\varphi\|$. It is easy to show that map $\left.\Phi\right|_{\widehat{\B}_0(\D)}$ is an isometric linear embedding from $\widehat{\B}_0(\D)$ into $\mathcal{C}_0(\D)$. By the Hahn--Banach Theorem, there exists a continuous linear functional $\widetilde{T}\colon \mathcal{C}_0(\D)\to\C$ such that $\widetilde{T}(\Phi(f))=T(\Phi(f))$ for all $f\in\widehat{\B}_0(\D)$ and $\|\widetilde{T}\|=\|T\|$. Now, by the Riesz Representation Theorem, there exists a finite regular Borel measure $\mu$ on $\D$ with total variation $\|\mu\|=\|\widetilde{T}\|$ such that
$$
\widetilde{T}(g)=\int_{\D}g\, d\mu
$$
for all $g\in\mathcal{C}_0(\D)$, and thus
$$
\varphi(f)=\int_{\D}\Phi(f)\, d\mu
$$
for all $f\in\widehat{\B}_0(\D)$. If we now define
$$
\gamma(f)=\int_{\D}\Phi(f)\, d\mu \qquad (f\in\widehat{\B}(\D)),
$$
it is clear that $\gamma\in\widehat{\B}(\D)^*$ and $\gamma(f)=\varphi(f)$ for all $f\in\widehat{\B}_0(\D)$. In light of Theorem \ref{teo-2}, Corollary \ref{cor-1} and the Banach--Dieudonn\'e Theorem, to ensure that $\gamma\in\G(\D)$, it suffices to prove that $\gamma$ is $\tau_p$-continuous on $B_{\widehat{\B}(\D)}$, where $\tau_p$ denotes the topology of pointwise convergence. Thus, let $(f_i)_{i\in I}$ be a net in $B_{\widehat{\B}(\D)}$ which converges pointwise on $\D$ to zero. Then $(\Phi(f_i))_{i\in I}$ converges pointwise on $\D$ to zero and, since $\left|\Phi(f_i)(z)\right|\leq\|\Phi(f_i)\|_{\infty}=\rho_\B(f_i)\leq 1$ for all $i\in I$ and for all $z\in\D$, we conclude that $(\gamma(f_i))_{i\in I}$ converges to $0$ by the Lebesgue's Bounded Convergence Theorem. Hence $\gamma\in\G(\D)$ and $R(\gamma)=\varphi$. 

Finally, fix $\gamma\in\G(\D)$ and let $f\in B_{\widehat{\B}(\D)}$. By Corollary \ref{cor-2}, there exists a net $(f_i)_{i\in I}$ in $B_{\widehat{\B}_0(\D)}$ which converges to $f$ in the topology $\tau_p$. Since $\gamma$ is $\tau_p$-continuous on $B_{\widehat{\B}(\D)}$ and  
$$
\left|\gamma(f_i)\right|=\left|R(\gamma)(f_i)\right|\leq\|R(\gamma)\|\rho_\B(f_i)\leq\|R(\gamma)\|
$$
for all $i\in I$, we have $\left|\gamma(f)\right|\leq\|R(\gamma)\|$ and so $\|\gamma\|\leq\|R(\gamma)\|$. This completes the proof.
\end{proof}

A Banach space $X$ is said to be a \textit{strongly unique predual of $X^*$} if for every Banach space $Y$, every isometric isomorphism from $X^*$ onto $Y^*$ is weak* continuous. Clearly, a strongly unique predual is a unique predual, which means that $X^*\cong Y^*$ implies $X\cong Y$. 

By \cite[III, Examples 1.4]{HarWerWer-93}, $\widehat{\B}_0(\D)$ is an $M$-ideal in $\widehat{\B}_0(\D)^{**}\cong\widehat{\B}(\D)$. In the light of Theorem \ref{messi1}, an application of \cite[III, Proposition 2.10]{HarWerWer-93} yields the following result (compare to \cite[Corollary 2]{Nar-90}).

\begin{proposition}\label{prop-unique}
$\G(\D)$ is a strongly unique predual of $\widehat{\B}(\D)$, and $\widehat{\B}_0(\D)$ is not a strongly unique predual of $\G(\D)$. $\hfill\Box$
\end{proposition}

Note that $B_{\widehat{\B}(\D)}$ is $w^*$-compact by Theorem \ref{teo-2}, and that, on $B_{\widehat{\B}(\D)}$, the weak* topology coincides with the topology of pointwise convergence $\tau_{p}$ by Corollary \ref{cor-1}. Thus the Dixmier--Ng Theorem \cite{Ng} assures that $\widehat{\B}(\D)$ is a dual Banach space with predual given by 
$$
\mathcal{G}_0(\D)=\{\gamma\in\widehat{\B}(\D)^*\colon \gamma|_{B_{\widehat{\B}(\D)}} \text{ is } \tau_p\text{-continuous}\}.
$$

Applying Proposition \ref{prop-unique}, we have the following identification.

\begin{corollary}
$\G(\D)$ is isometrically isomorphic to $\mathcal{G}_0(\D)$.$\hfill\Box$
\end{corollary}


\section{Linearization of derivatives of Banach-valued Bloch mappings}\label{section 3}

Next we state the universal extension property of $\G(\D)$.

\begin{theorem}\label{teo-3}
For every complex Banach space $X$ and every Bloch mapping $f\colon\D\to X$ with $f(0)=0$, there exists a unique bounded linear operator $S_f\colon\G(\D)\to X$ such that $S_f\circ\Gamma=f'$, that is, the diagram 
$$
\begin{tikzpicture}
\node (D) {$\D$};
\node (GD) [below of=D] {$\G(\D)$};
\node (X) [right of=GD] {$X$};
\draw[->] (D) to node {$f'$} (X);
\draw[->] (D) to node [swap] {$\Gamma$} (GD);
\draw[->, dashed] (GD) to node [swap] {$S_f$} (X);
\end{tikzpicture}
$$
commutes. Furthermore, $\|S_f\|=\rho_{\B}(f)$.  
\end{theorem}

\begin{proof}
Let $X$ be a complex Banach space and let $f\in\widehat{\B}(\D,X)$. Define the mapping $R_f\colon\lin(\Gamma(\D))\to X$ by 
$$
R_f(\gamma)=\sum_{k=1}^n\lambda_kf'(z_k)\qquad \left(\gamma=\sum_{k=1}^n\lambda_k\gamma_{z_k}\right).
$$
Clearly, $R_f$ is linear and $R_f\circ\Gamma=f'$. Let $\gamma\in\lin(\Gamma(\D))$. For every $x^*\in B_{X^*}$, we obtain
$$
|x^*(R_f(\gamma))| = |\gamma(x^*\circ f)| \leq \|\gamma\|\rho_{\B}(x^*\circ f) \leq \|\gamma\|\rho_{\B}(f).
$$
Taking the supremum over all $x^*\in B_{X^*}$ yields
$$
\|R_f(\gamma)\| \leq\rho_{\B}(f)\|\gamma\|.
$$
Hence $R_f$ is continuous with $\|R_f\|\leq \rho_{\B}(f)$. Since $\lin(\Gamma(\D))$ is norm-dense in $\G(\D)$, there exists a unique continuous mapping $S_f\colon\G(\D)\to X$ such that $S_f(\gamma)=R_f(\gamma)$ for all $\gamma\in\lin(\Gamma(\D))$. Further, $S_f$ is linear and $\|S_f\|=\|R_f\|$. Hence $S_f\circ\Gamma=f'$ and $\|S_f\|\leq \rho_{\B}(f)$. Conversely, we have $\rho_{\B}(f)\leq\|S_f\|$ because  
$$
\|S_f\|\geq\|S_f(\widehat{\gamma}_z)\|=(1-|z|^2)\|S_f(\gamma_z)\|=(1-|z|^2)\|(S_f\circ\Gamma)(z)\|=(1-|z|^2)\|f'(z)\|
$$
for all $z\in\D$. This completes the proof of the theorem. 
\end{proof}

Theorem \ref{teo-3} allows to identify the spaces $\widehat{\B}(\D,X)$ and $\L(\G(\D),X)$. 

\begin{corollary}\label{cor-3-new}
The mapping $f\mapsto S_f$ is an isometric isomorphism from $\widehat{\B}(\D,X)$ onto $\L(\G(\D),X)$.
\end{corollary}

\begin{proof}
Let $f,g\in\widehat{\B}(\D,X)$ and $\alpha,\beta\in\C$. Using Theorem \ref{teo-3}, we have 
$$
(\alpha S_f+\beta S_g)\circ\Gamma=\alpha(S_f\circ\Gamma)+\beta(S_g\circ\Gamma)=\alpha f'+\beta g'=(\alpha f+\beta g)'
$$
and the cited uniqueness in that theorem yields $S_{\alpha f+\beta g}=\alpha S_f+\beta S_g$. Hence the mapping of the statement is linear, and it is also isometric by Theorem \ref{teo-3}. To prove its surjectivity, take $T\in\L(\G(\D),X)$. Clearly, $T\circ\Gamma\in\H(\D,X)$, and Lemma \ref{main lemma} gives us a holomorphic function $f\colon\D\to X$ with $f(0)=0$ such that $f'=T\circ\Gamma$. Moreover, $f$ is Bloch since 
$$
(1-|z|^2)|f'(z)|=(1-|z|^2)\|T(\gamma_z)\|=\|T(\widehat{\gamma}_z)\|\leq \|T\|
$$
for all $z\in\D$. We have
$$
T(\gamma_z)=T\circ\Gamma(z)=f'(z)=S_f(\gamma_z)
$$
for all $z\in\D$. By 
the norm-denseness of $\lin(\Gamma(\D))$ in $\G(\D)$, we conclude that $T=S_f$. 
\end{proof}

The space $\H^\infty(\D,X)$ of all bounded holomorphic mappings $f\colon\D\to X$ is contained in the space of all Bloch mappings from $\D$ into $X$ with $\rho_{\B}(f)\leq \|f\|_\infty$, where 
$$
\|f\|_\infty=\sup\left\{\|f(z)\|\colon z\in\D\right\},
$$
but there are unbounded Bloch functions such as the function $z\mapsto\log(1-z)$ from $\D$ to $\C$. 

Mujica \cite{Muj-91} constructed a Banach space $\G^\infty(\D)$ and a mapping $g_\D\in\H^\infty(\D,\G^\infty(\D))$ with the following universal property: for each complex Banach space $X$ and each mapping $f\in\H^\infty(\D,X)$, there exists a unique operator $T_f\in\L(\G^\infty(\D),X)$ such that $T_f\circ g_\D=f$ and $\|T_f\|=\|f\|_\infty$. 

The space $\G^\infty(\D)$ is a predual of $\H^\infty(\D)$ and it is defined in \cite{Muj-91} as the closed subspace of all functionals $\phi\in\H^\infty(\D)^*$ such that the restriction of $\phi$ to $B_{\H^\infty(\D)}$ is $\tau_c$-continuous, where $\tau_c$ is the compact-open topology. The mapping $g_\D\colon\D\to\G^\infty(\D)$ is defined by $g_\D(z)=\delta_z$ for $z\in\D$, where $\delta_z$ is the evaluation functional at $z$ defined on $\H^\infty(\D)$.

Next we show the connection between our linearization $S_f\in\L(\G(\D),X)$ and the Mujica's linearization $T_f\in\L(\G^\infty(\D),X)$ of a zero-preserving mapping $f\in\H^\infty(\D,X)$.

\begin{proposition}
For every $f\in\H^\infty(\D,X)$ with $f(0)=0$, there are operators $S_f\in\L(\G(\D),X)$ and $T_f\in\L(\G^\infty(\D),X)$ with $\|S_f\|\leq \|T_f\|$ such that $S_f\circ\Gamma=(T_f\circ g_\D)'$.
\end{proposition}

\begin{proof}
Let $f\in\H^\infty(\D,X)$ with $f(0)=0$. On a hand, there exists $T_f\in\L(\G^\infty(\D),X)$ with $\|T_f\|=\|f\|_\infty$ such that $T_f\circ g_\D=f$ by \cite[Theorem 2.1]{Muj-91}. On the other hand, since $f\in\widehat{\B}(\D,X)$, Theorem \ref{teo-3} provides $S_f\in\L(\G(\D),X)$ with $\|S_f\|=\rho_\B(f)$ such that $S_f\circ\Gamma=f'$. Hence $S_f\circ\Gamma=(T_f\circ g_\D)'$, and since $\rho_\B(f)\leq \|f\|_\infty$, we have $\|S_f\|\leq \|T_f\|$. 
\end{proof}

We may study holomorphic mappings on Banach spaces in terms of power series expansions (see \cite[Appendix A]{AreBatHieNeu-01}). 

Let $X$ be a Banach space and let $\Omega\subseteq\C$ be an open set. If $f\colon\Omega\to X$ is holomorphic, by Taylor's Theorem for vector-valued holomorphic mappings,  
for each $a\in\Omega$ there exist a closed disc $\overline{D}(a,r_a)\subseteq\Omega$ and a sequence $\{a_m^{(a)}\}_{m\geq 0}$ in $X$ given by 
$$
a_m^{(a)}=\frac{1}{2\pi i}\int_{|z-a|=r_a}\frac{f(z)}{(z-a)^{m+1}}\ dz,
$$
such that  
$$
f(z)=\sum_{m=0}^\infty a_m^{(a)}(z-a)^m,
$$
where the series converges uniformly for $z\in D(a,r_a)$. This series is known as the \textit{Taylor series of $f$ at $a$} and by analogy with the complex case, it is represented as 
$$
f(z)=\sum_{m=0}^\infty \frac{1}{m!}d^mf(a)(z-a).
$$

\begin{proposition}
$\C$ is topologically isomorphic to a complemented subspace of $\G(\D)$.
\end{proposition}

\begin{proof}
Let $f\colon\D\to\C$ be the polynomial defined by $f(z)=z^2/2$ for all $z\in\D$. Clearly, $f\in\widehat{\B}(\D)$ with $\rho_\B(f)\leq 1$. 
By Theorem \ref{teo-3}, there exists $S_f\in\G(\D)^*$ with $\|S_f\|=\rho_\B(f)$ such that 
$$
S_f\circ\Gamma(z)=f'(z)=z\qquad (z\in\D).
$$
Since $\Gamma\colon\D\to\G(\D)$ is holomorphic, let $S=d^1\Gamma(0)\in\L(\C,\G(\D))$. By Cauchy's integral formula for derivatives, 
we have 
$$
S(z)=\frac{1}{2\pi i}\int_{|w|=r}\frac{\Gamma(wz)}{w^2}\ dw\qquad (z\in\C),
$$
where $r$ is chosen so that $\left\{wz\colon |w|\leq r\right\}\subseteq\D$. It follows that 
\begin{align*}
S_f\circ S(z)&=S_f\circ d^1\Gamma(0)(z)=d^1(S_f\circ\Gamma)(0)(z)\\
             &=\frac{1}{2\pi i}\int_{|w|=r}\frac{S_f\circ\Gamma(wz)}{w^2}\ dw=\frac{1}{2\pi i}\int_{|w|=r}\frac{wz}{w^2}\ dw\\
						 &=\frac{z}{2\pi i}\int_{|w|=r}\frac{dw}{w}=z \qquad (z\in\C).
\end{align*}
Hence 
$$
|z|=|S_f\circ S(z)|\leq \|S_f\|\|S(z)\|\leq\|S(z)\|
\qquad (z\in\C),
$$
and so $S$ is injective and $S^{-1}\colon S(\C)\to\C$ is continuous. 
Define $P=S\circ S_f$. Clearly, $P$ is a continuous linear projection from $\G(\D)$ into itself. 
Moreover, $P(\G(\D))=S(S_f(\G(\D)))=S(\C)$. Hence $S$ is a topological isomorphism from $\C$ onto the complemented subspace $P(\G(\D))$ of $\G(\D)$. 

\end{proof}


\section{Bloch mappings with compact type Bloch range}\label{section 4}

A linear operator between Banach spaces $T\colon X\to Y$ is said to be \textit{compact (weakly compact)} if $T(B_X)$ is relatively compact (respectively, relatively weakly compact) in $Y$. We denote by $\K(X,Y)$, $\W(X,Y)$, $\F(X,Y)$ and $\overline{\F}(X,Y)$ the spaces of compact linear operators, weakly compact linear operators, bounded finite-rank linear operators and approximable linear operators from $X$ into $Y$, respectively. 

In this section, we introduce and study the analogues of these concepts in the setting of Bloch spaces.

\begin{definition}\label{def-3-1}
Let $X$ be a complex Banach space and $f\in\H(\D,X)$. The \textit{Bloch range} of $f$, denoted $\rang_{\B}(f)$, is the set
$$
\left\{(1-|z|^2)f'(z)\colon z\in\D\right\}\subseteq X.
$$
Note that $f$ is Bloch if and only if $\rang_{\B}(f)$ is a bounded subset of $X$. 

The mapping $f$ is said to be \textit{compact Bloch (weakly compact Bloch)} if $\rang_{\B}(f)$ is a relatively compact (resp., relatively weakly compact) subset of $X$.

We denote by $\widehat{\B}_{\K}(\D,X)$ and $\widehat{\B}_{\W}(\D,X)$ the linear spaces of compact Bloch mappings and weakly compact Bloch mappings from $\D$ to $X$ that preserve the zero, respectively.
\end{definition}

Other classes of Bloch mappings closely related to the previous ones are the following.

\begin{definition}
Let $X$ be a complex Banach space. A Bloch mapping $f\in\widehat{\B}(\D,X)$ has \textit{finite dimensional Bloch rank} if the linear hull of the set $\left\{(1-|z|^2)f'(z)\colon z\in\D\right\}$ is a finite dimensional subspace of $X$. In that case, we define the \textit{Bloch rank $\B$-$\rank(f)$ of $f$} to be the dimension of this subspace. 

A Bloch mapping $f\in\widehat{\B}(\D,X)$ is said to be \textit{approximable} if it is the limit in the Bloch norm $\rho_{\B}$ of a sequence of finite-rank Bloch mappings of $\widehat{\B}(\D,X)$.

We denote by $\widehat{\B}_{\F}(\D,X)$ and $\widehat{\B}_{\overline{\F}}(\D,X)$ the spaces of finite-rank Bloch mappings and approximable Bloch mappings $f$ from $\D$ to $X$ such that $f(0)=0$, respectively.
\end{definition}


\subsection{Invariance by M\"obius transformations}

The \textit{M\"obius group of $\D$}, denoted by $\Aut(\D)$, consists of all one-to-one holomorphic functions $\phi$ that map $\D$ onto itself. Each $\phi\in\Aut(\D)$ has the form $\phi=\lambda\phi_a$ with $\lambda\in\T$ and $a\in\D$, where
$$
\phi_a(z)=\frac{a-z}{1-\overline{a}z}\qquad (z\in\D).
$$
It is easy to check that  
$$
\left|\phi_a'(z)\right|=\frac{1-|a|^2}{|1-\overline{a}z|^2}
$$
and
$$
1-|\phi_a(z)|^2=\frac{(1-|a|^2)(1-|z|^2)}{|1-\overline{a}z|^2}=(1-|z|^2)\left|\phi_a'(z)\right|
$$
for all $z\in\D$. 

Let us recall (see \cite{AraFisPee-85}) that for a complex Banach space $X$, a linear space $\A(\D,X)$ of holomorphic mappings from $\D$ into $X$ endowed with a seminorm $\rho_\A$ is said to be \textit{M\"obius-invariant} if 
\begin{enumerate}
	\item $\A(\D,X)\subseteq \B(\D,X)$ and there exists a constant $c>0$ such that $\rho_\B(f)\leq c\rho_\A(f)$ for all $f\in\A(\D,X)$,
	\item $f\circ\phi\in\A(\D,X)$ with $\rho_\A(f\circ\phi)=\rho_\A(f)$ for all $\phi\in\Aut(\D)$ and $f\in\A(\D,X)$.
\end{enumerate}

For $f\in\H(\D,X)$ and $\phi\in\Aut(\D)$, we have  
$$
(1-|z|^2)(f\circ\phi)'(z)=(1-|z|^2)f'(\phi(z))\phi'(z)=(1-\left|\phi(z)\right|^2)f'(\phi(z))\frac{\phi'(z)}{|\phi'(z)|}
$$
for all $z\in\D$. Hence $\rang_{\B}(f\circ\phi)\subseteq\T\, \rang_{\B}(f)$, and therefore, if $f$ is (weakly) compact Bloch, it is clear that $f\circ\phi$ is (weakly) compact Bloch with $\rho_\B(f\circ\phi)=\rho_\B(f)$. Hence we have proved the following.

\begin{proposition}\label{prop-3-1}
The spaces of compact Bloch mappings and weakly compact Bloch mappings from $\D$ to $X$ are M\"obius-invariant. $\hfill\qed$
\end{proposition}


\subsection{Linearization and factorization}

We now study the relationship between the Bloch compactness of $f\in\widehat{\B}(\D,X)$ and the compactness of its linearization $S_f\in\L(\G(\D),X)$. 

\begin{theorem}\label{teo-3-1}
Let $f\in\widehat{\B}(\D,X)$. The following conditions are equivalent:
\begin{enumerate}
\item $f$ is compact Bloch.
\item $S_f\colon\G(\D)\to X$ is compact.
\end{enumerate}
Furthermore, the mapping $f\mapsto S_{f}$ is an isometric isomorphism from $\widehat{\B}_{\K}(\D,X)$ onto $\K(\G(\D),X)$.
\end{theorem}

\begin{proof}
First, using Theorem \ref{teo-3}, we derive the relations:
\begin{align*}
\rang_{\B}(f)=S_f(\M_\B(\D))&\subseteq S_f(B_{\G(\D)})=S_f(\overline{\abco}(\M_\B(\D)))\\
						 &\subseteq \overline{\abco}(S_f(\M_\B(\D)))=\overline{\abco}(\rang_{\B}(f))\\
						 &\subseteq\overline{\abco}(\overline{\rang_{\B}(f)})
\end{align*}

$(1)\Rightarrow (2)$: If $f\in\widehat{\B}_{\K}(\D,X)$, we deduce that $S_f\in\K(\G(\D),X)$ by applying the second inclusion above and the fact that $\overline{\abco}(\overline{\rang_{\B}(f)})$ is compact in $X$ by Mazur's Compactness Theorem. 

$(2)\Rightarrow (1)$: If $S_f\in\K(\G(\D),X)$, we obtain that $f\in\widehat{\B}_{\K}(\D,X)$ by using the first inclusion above. 

The last assertion of the statement follows immediately using Corollary \ref{cor-3-new} and what was proved above.
\end{proof}

As an application of Theorem \ref{teo-3-1}, we show that the members of the Bloch space $\widehat{\B}_{\K}$ are generated by composition with the Banach operator ideal $\K$.

\begin{corollary}\label{messi-3}
Let $f\in\widehat{\B}(\D,X)$. The following conditions are equivalent:
\begin{enumerate}
\item $f\colon \D\to X$ is compact Bloch. 
\item $f=T\circ g$, where $Y$ is a complex Banach space, $g\in\widehat{\B}(\D,Y)$ and $T\in\K(Y,X)$. 
\end{enumerate}
In this case, $\rho_{\B}(f)=\inf\{\|T\|\rho_\B(g)\}$, where the infimum is taken over all factorizations of $f$ as above, and this infimum is attained at $T=S_f$ and $g=\Gamma$.
\end{corollary}

\begin{proof}
$(1)\Rightarrow (2)$: If $f\in\widehat{\B}_{\K}(\D,X)$, we have the factorization $f'=S_f\circ\Gamma$, where $\G(\D)$ is a complex Banach space, $S_f\in\K(\G(\D),X)$ and $\Gamma\in\H(\D,\G(\D))$ by Theorems \ref{teo-3} and \ref{teo-3-1}. Now, by Lemma \ref{main lemma}, there is a holomorphic mapping $h\colon\D\to\G(\D)$ with $h(0)=0$ satisfying that $h'(z)=\Gamma(z)$ for all $z\in\D$. Since 
$$
(1-|z|^2)\|h'(z)\|=(1-|z|^2)\|\Gamma(z)\|=1
$$
for all $z\in\D$, we get that $h\in\widehat{\B}(\D,\G(\D))$ with $\rho_{\B}(h)=1$. Thus we obtain the factorization $f'=S_f\circ h'=(S_f\circ h)'$. This implies $f=S_f\circ h$ since $\D$ is connected and $f(0)=(S_f\circ h)(0)=0$. Moreover, we have $\inf\left\{\|T\|\rho_\B(g)\right\}\leq \|S_f\|\rho_\B(h)=\rho_{\B}(f)$.

$(2)\Rightarrow (1)$: Assume that $f=T\circ g$, where $Y$ is a complex Banach space, $g\in\widehat{\B}(\D,Y)$ and $T\in\K(Y,X)$. Then $f'=T\circ g'$. Since $g'=S_g\circ\Gamma$ by Theorem \ref{teo-3}, it follows that $f'=T\circ S_g\circ\Gamma$ which implies that $S_f=T\circ S_g$, and thus $S_f\in\K(\G(\D),X)$ by the ideal property of $\K$. By Theorem \ref{teo-3-1}, we obtain that $f\in\widehat{\B}_{\K}(\D,X)$ with  
$$
\rho_{\B}(f)=\|S_f\|=\|T\circ S_g\|\leq \|T\|\|S_g\|=\|T\|\rho_\B(g), 
$$
and taking the infimum over all representations of $f$, we deduce that $\rho_{\B}(f)\leq\inf\{\|T\|\rho_\B(g)\}$. 
\end{proof}

The Davis--Figiel--Johnson--Pe\l czy\'nski Factorization Theorem \cite{dfjp} states that any weakly compact linear operator between Banach spaces factors through a reflexive Banach space. We now extend this result to weakly compact Bloch mappings.

\begin{theorem}\label{teo-3-2}
Let $f\in\widehat{\B}(\D,X)$. The following are equivalent:
\begin{enumerate}
	\item $f\colon\D\to X$ is Bloch weakly compact.
	\item $S_f\colon\G(\D)\to X$ is weakly compact.
	\item There exist a reflexive Banach space $Y$, an operator $T\in\L(Y,X)$ and a mapping $g\in\widehat{\B}(\D,Y)$ such that $f=T\circ g$.
\end{enumerate}
Furthermore, $f\mapsto S_f$ is an isometric isomorphism from $\widehat{\B}_\W(\D,X)$ onto $\W(\G(\D),X)$.
\end{theorem}

\begin{proof}
Since the norm closure and weak closure of the convex hull of a subset of a normed space coincide \cite[Corollary 2.5.18]{Meg-98} and the norm-closed convex hull of a weakly compact subset of a Banach space is itself weakly compact \cite[Theorem 2.8.14]{Meg-98}, a similar proof to that of Theorem \ref{teo-3-1} yields that $(1)\Leftrightarrow(2)$.  

If (2) holds, applying the Davis--Figiel--Johnson--Pe\l czy\'nski Theorem, there exist a reflexive Banach space $Y$ and operators $T\in\L(Y,X)$ and $S\in\L(\G(\D),Y)$ such that $S_f=T\circ S$. Since $S\circ\Gamma\in\H(\D,Y)$, by Lemma \ref{main lemma} we can find a mapping $g\in\H(\D,Y)$ with $g(0)=0$ such that $g'=S\circ\Gamma$. It follows that $g\in\widehat{\B}(\D,Y)$ since 
$$
(1-|z|^2)\|g'(z)\|=(1-|z|^2)\|S(\Gamma(z))\|\leq (1-|z|^2)\|S\|\|\Gamma(z)\|=\|S\|
$$
for all $z\in\D$. Moreover, $f'=S_f\circ\Gamma=T\circ S\circ\Gamma=T\circ g'=(T\circ g)'$, hence $f=T\circ g$ and this proves (3). 

Finally, (3) implies (1) because $\rang_\B(f)=T(\rang_\B(g))$,  
where $T$ is weak-to-weak continuous 
and $\rang_\B(g)$ 
is relatively weakly compact in $Y$ since it is a bounded subset of the reflexive Banach space $Y$. 
\end{proof} 

Let us recall that if $A$ is a set and $X$ is a linear space, then a mapping $f\colon A\to X$ is said to have \textit{finite dimensional rank} if the linear hull of its range is a finite dimensional subspace of $X$ in whose case the \textit{rank of $f$}, denoted by $\rank(f)$, is defined as the dimension of $\lin(f(A))$. 

\begin{theorem}\label{prop2.1}
Let $f\in\widehat{\B}(\D,X)$. The following are equivalent:
\begin{enumerate}
  \item $f$ has finite dimensional Bloch rank.
	\item $f'$ has finite dimensional rank.
	\item $S_f\in\L(\G(\D),X)$ has finite rank.
\end{enumerate}
In that case, $\lin(f'(\D))=S_f(\G(\D))$ and $\B$-$\rank(f)=\rank(f')=\rank(S_f)$. Furthermore, the mapping $f\mapsto S_f$ is an isometric isomorphism from $\widehat{\B}_\F(\D,X)$ onto $\F(\G(\D),X)$.
\end{theorem}

\begin{proof}
First, observe that 
$$
\lin\left\{(1-|z|^2)f'(z)\colon z\in\D\right\}=\lin\left\{f'(z)\colon z\in\D\right\}
$$	
for any function $f\in\widehat{\B}(\D,X)$. 

$(1)\Leftrightarrow(2)$ and $\B$-$\rank(f)=\rank(f')$ follow from this observation. We now prove $(2)\Leftrightarrow(3)$. If $f'$ has finite dimensional rank, then $\lin(f'(\D))$ is finite dimensional, hence it is closed in $X$. 
Using Theorem \ref{teo-3}, we have
$$
S_f(\G(\D))=S_f(\overline{\lin}(\Gamma(\D)))\subseteq\overline{S_f(\lin(\Gamma(\D)))}=\overline{\lin}(S_f(\Gamma(\D)))=\overline{\lin}(f'(\D))=\lin(f'(\D))
$$
and hence $S_f$ has finite rank. Conversely, if $S_f$ has finite rank, then $f'$ has finite dimensional rank since
$$
\lin(f'(\D))=\lin(S_f(\Gamma(\D)))=S_f(\lin(\Gamma(\D)))\subseteq S_f(\overline{\lin}(\Gamma(\D)))=S_f(\G(\D)).
$$
\end{proof}

\subsection{Inclusions}

Note that $\widehat{\B}_{\F}(\D,X)$ is a linear subspace of $\widehat{\B}_{\K}(\D,X)$. In fact, we have: 

\begin{proposition}
Every approximable Bloch mapping from $\D$ to $X$ is compact Bloch.
\end{proposition}

\begin{proof}
Let $f\in\widehat{\B}_{\overline{\F}}(\D,X)$. Then there is a $(f_n)_{n\in\N}$ in $\widehat{\B}_{\F}(\D,X)$ such that $\rho_\B(f_n-f)\to 0$ as $n\to\infty$. Since $S_{f_n}\in\F(\G(\D),X)$ by Theorem \ref{prop2.1}, $\F(\G(\D),X)\subseteq\K(\G(\D),X)$ and $\|S_{f_n}-S_f\|=\|S_{f_n-f}\|=\rho_\B(f_n-f)$ for all $n\in\N$ by Corollary \ref{cor-3-new}
, we deduce that $S_f\in\K(\G(\D),X)$ by \cite[Corollary 3.4.9]{Meg-98}, and so $f\in\widehat{\B}_{\K}(\D,X)$ by Theorem \ref{teo-3-1}.
\end{proof}

An application of Theorem \ref{prop2.1} and Corollary \ref{cor-3-new} shows the connection of an approximable Bloch mapping with its linearization:

\begin{theorem}\label{prop2.1.new}
Let $f\in\widehat{\B}(\D,X)$. Then $f$ is approximable Bloch if and only if $S_f\in\L(\G(\D),X)$ is approximable. Furthermore, the mapping $f\mapsto S_f$ is an isometric isomorphism from $\widehat{\B}_{\overline{\F}}(\D,X)$ onto $\overline{\F}(\G(\D),X)$. $\hfill\qed$
\end{theorem}

Let us recall that a Banach space $X$ is said to have the \textit{approximation property} if given a compact set $K\subseteq X$ and $\varepsilon>0$, there is an operator $T\in\F(X,X)$ such that $\|T(x)-x\|<\varepsilon$ for every $x\in K$. On the approximation property in holomorphic function spaces, we refer to results of Aron and Schottenloher \cite{AroSch-76} and Mujica \cite{Muj-91}.

Grothendieck \cite{g} proved that a dual Banach space $X^*$ has the approximation property if and only if given a Banach space $Y$, an operator $S\in\K(X,Y)$ and $\varepsilon>0$, there is an operator $T\in\F(X,Y)$ such that $\|T-S\|<\varepsilon$. Since $\widehat{\B}(\D)$ is a dual space, combining the Grothendieck's result with Theorems \ref{prop2.1} and \ref{teo-3-1}, we next give a necessary and sufficient condition for the space $\widehat{\B}(\D)$ have the approximation property.

\begin{corollary}
$\widehat{\B}(\D)$ has the approximation property if and only if, for each complex Banach space $X$, every compact Bloch mapping from $\D$ to $X$ is approximable Bloch. $\hfill\qed$
\end{corollary}


\subsection{Banach ideal property}

We now formalize the notion of an ideal of Bloch mappings with a definition inspired by the notion of operator ideal between Banach spaces \cite{Pie-80}.

\begin{definition}\label{def-ideal}

An \textit{ideal of normalized Bloch mappings} (or simply, a \textit{normalized Bloch ideal}) is a subclass $\I^{\widehat{\B}}$ of the class of all normalized Bloch mappings $\widehat{\B}$ such that for each complex Banach space $X$, the components 
$$
\I^{\widehat{\B}}(\D,X)=\I^{\widehat{\B}}\cap\widehat{\B}(\D,X)
$$ 
satisfy:
\begin{enumerate}
\item[(I1)] $\I^{\widehat{\B}}(\D,X)$ is a linear subspace of $\widehat{\B}(\D,X)$,
\item[(I2)] For any $g\in\widehat{\B}(\D)$ and $x\in X$, the mapping $g\cdot x\colon z\mapsto g(z)x$ from $\D$ to $X$ is in $\I^{\widehat{\B}}(\D,X)$,
\item[(I3)] \textit{The ideal property}: if $f\in\I^{\widehat{\B}}(\D,X)$, $h\colon\D\to\D$ is a holomorphic function with $h(0)=0$ and $T\in\L(X,Y)$ where $Y$ is a complex Banach space, then $T\circ f\circ h$ belongs to $\I^{\widehat{\B}}(\D,Y)$.
\end{enumerate}
A normalized Bloch ideal $\I^{\widehat{\B}}$ is said to be \textit{normed (Banach)} if there is a function $\|\cdot\|_{\I^{\widehat{\B}}}\colon\I^{\widehat{\B}}\to\mathbb{R}_0^+$ such that for every complex Banach space $X$, the following three conditions are satisfied:
\begin{enumerate}
	\item[(N1)] $(\I^{\widehat{\B}}(\D,X),\|\cdot\|_{\I^{\widehat{\B}}})$ is a normed (Banach) space with $\rho_\B(f)\leq\|f\|_{\I^{\widehat{\B}}}$ for all $f\in\I^{\widehat{\B}}(\D,X)$,
	\item[(N2)] $\|g\cdot x\|_{\I^{\widehat{\B}}}=\rho_\B(g)\|x\|$ for all $g\in\widehat{\B}(\D)$ and $x\in X$, 
	\item[(N3)] If $Y$ is a complex Banach space, $f\in\I^{\widehat{\B}}(\D,X)$, $h\colon\D\to\D$ is a holomorphic function with $h(0)=0$ and $T\in\L(X,Y)$, then $\|T\circ f \circ h\|_{\I^{\widehat{\B}}}\leq \|T\|\|f\|_{\I^{\widehat{\B}}}$.
\end{enumerate}
A normed normalized Bloch ideal $[\I^{\widehat{\B}},\|\cdot\|_{\I^{\widehat{\B}}}]$ is said to be:
\begin{enumerate}
\item[(I)] \textit{Injective} if for any mapping $f\in\widehat{\B}(\D,X)$, any complex Banach space $Y$ and any metric injection $\iota\in\L(X,Y)$, we have that $f\in\I^{\widehat{\B}}(\D,X)$ with $\|f\|_{\I^{\widehat{\B}}}=\|\iota\circ f\|_{\I^{\widehat{\B}}}$ whenever $\iota\circ f\in\I^{\widehat{\B}}(\D,Y)$.
\item[(S)] \textit{Surjective} if for any complex Banach space $X$, any $f\in\widehat{\B}(\D,X)$ and any holomorphic function $\pi\colon\D\to\D$ with $\pi(0)=0$ such that $\widehat{\pi}\in\L(\G(\D),\G(\D))$ is a metric surjection, we have that $f\in\I^{\widehat{\B}}(\D,X)$ with $\|f\|_{\I^{\widehat{\B}}}=\|f\circ \pi\|_{\I^{\widehat{\B}}}$ whenever $f\circ \pi\in\I^{\widehat{\B}}(\D,X)$.
\end{enumerate}
\end{definition}


The condition imposed on $\widehat{\pi}$ in the surjectivity property is closely related to isometric composition operators on $\widehat{\B}(\D)$, studied by Mart\'in and Vukoti\'c \cite{MarVuk-07}. 

\begin{remark}
In view of Corollary \ref{cor-august} and \cite[Theorem 1.1]{MarVuk-07}, if $h\colon\D\to\D$ is a holomorphic function such that $h(0)=0$, then the following are equivalent:
\begin{enumerate}
	\item $\widehat{h}\colon\G(\D)\to\G(\D)$ is a metric surjection.
	\item $C_h\colon\widehat{\B}(\D)\to\widehat{\B}(\D)$ is a metric injection. 
	\item Either $h$ is a rotation of the identity or $h$ has the property (M): for every $a\in\D$ there exists a sequence $(z_n)$ in $\D$ such that $(|z_n|)\to 1$, $(h(z_n))\to a$, and $((1-|z_n|^2)|h'(z_n)|/(1-|h(z_n)|^2))\to 1$.
\end{enumerate}
\end{remark}

We present some examples of such normalized Bloch ideals.

\begin{proposition}\label{ideal-0}
$[\widehat{\B},\rho_\B]$ is a Banach ideal of normalized Bloch mappings. Furthermore, if $X$ and $Y$ are complex Banach spaces, $T\in\L(X,Y)$, $h\colon\D\to\D$ is a holomorphic function with $h(0)=0$ and $f\in\widehat{\B}(\D,X)$, then $S_{T\circ f\circ h}=T\circ S_f\circ\widehat{h}$.
\end{proposition}

\begin{proof}
By Corollary \ref{cor-3-new}, $(\widehat{\B}(\D,X),\rho_\B)$ is a Banach space. Given $g\in\widehat{\B}(\D)$ and $x\in X$, it is clear that $g\cdot x\colon\D\to X$ is holomorphic with $(g\cdot x)'=g'\cdot x$ and $(g\cdot x)(0)=g(0)x=0$. Moreover, $g\cdot x$ is Bloch with $\rho_\B(g\cdot x)=\rho_\B(g)\|x\|$. Now, let $T$, $h$ and $f$ be as in the statement. Clearly, $(T\circ f\circ h)(0)=0$ and $T\circ f\circ h\in\H(\D,Y)$ with 
$$
(T\circ f\circ h)'=T\circ(f\circ h)'=T\circ[h'\cdot (f'\circ h)],
$$
and since 
$$
(1-|z|^2)|h'(z)|\leq 1-|h(z)|^2\qquad (z\in\D)
$$
by the Pick--Schwarz Lemma, we deduce that $T\circ f\circ h\in\widehat{\B}(\D,Y)$ with $\rho_\B(T\circ f\circ h)\leq \|T\|\rho_\B(f)$. This completes the proof of the first assertion of the statement. To prove the second assertion, note that, by Corollary \ref{cor-3}, there is an operator $\widehat{h}\in\L(\G(\D),\G(\D))$ with $\|\widehat{h}\|\leq 1$ such that $\widehat{h}\circ\Gamma=h'\cdot(\Gamma\circ h)$.  Since  
\begin{align*}
(T\circ f\circ h)'&=T\circ[h'\cdot (f'\circ h)]=T\circ[h'\cdot(S_f\circ\Gamma\circ h)]\\
                  &=T\circ [S_f(h'\cdot(\Gamma\circ h))]=T\circ[S_f\circ(\widehat{h}\circ\Gamma)]\\
									&=(T\circ S_f\circ\widehat{h})\circ\Gamma
\end{align*}
with $T\circ S_f\circ\widehat{h}\in\L(\G(\D),Y)$, we get that $S_{T\circ f\circ h}=T\circ S_f\circ\widehat{h}$ by Theorem \ref{teo-3}. 
\end{proof}

\begin{proposition}\label{ideal}
$[\widehat{\B}_\I,\rho_\B]$ is a Banach ideal of normalized Bloch mappings whenever $\I=\K,\W,\overline{\F}$. Furthermore, the ideal $[\widehat{\B}_\I,\rho_\B]$ is surjective and injective for $\I=\K,\W$.
\end{proposition}

\begin{proof}
Let $X$ be a complex Banach space and $\I=\K,\W,\overline{\F}$. We will show that $(\widehat{\B}_{\I}(\D,X),\rho_{\B})$ enjoys the required conditions: 

(N1) $(\widehat{\B}_\I(\D,X),\rho_\B)$ is a Banach space by Theorems \ref{teo-3-1}, \ref{teo-3-2} and \ref{prop2.1.new}. 

(N2) Given $g\in\widehat{\B}(\D)$ and $x\in X$, it is clear that $\rang_\B(g\cdot x)=\rang_\B(g)x$ and thus $g\cdot x\in\widehat{\B}_\F(\D,X)\subseteq\widehat{\B}_\I(\D,X)$ with $\rho_\B(g\cdot x)=\rho_\B(g)\|x\|$. 

(N3) Let $T\in\L(X,Y)$, $h\colon\D\to\D$ be a holomorphic function with $h(0)=0$ and $f\in\widehat{\B}_\I(\D,X)$. Then $S_f\in\I(\G(\D),X)$ by Theorems \ref{teo-3-1}, \ref{teo-3-2} and \ref{prop2.1.new}, hence $S_{T\circ f\circ h}=T\circ S_f\circ\widehat{h}\in\I(\G(\D),Y)$ by Proposition \ref{ideal-0} and the operator ideal property of $\I$, and thus $T\circ f\circ h\in\widehat{\B}_\I(\D,Y)$ with 
$$
\rho_\B(T\circ f\circ h)=\|S_{T\circ f\circ h}\|=\|T\circ S_f\circ\widehat{h}\|\leq\|T\|\rho_\B(f)
$$
again by the aforementioned theorems. 

Assume now that $\I=\K,\W$.

(I) Let $f\in\widehat{\B}(\D,X)$ and let $\iota\in\L(X,Y)$ be a metric injection, where $Y$ is a complex Banach space. Assume that $\iota\circ f\in\widehat{\B}_\I(\D,Y)$. Then $\iota\circ S_f=S_{\iota\circ f}\in\I(\G(\D),Y)$. Since the operator ideal $[\I,\|\cdot\|]$ is injective by \cite[4.6.12]{Pie-80}, it follows that $S_f\in\I(\G(\D),X)$ with $\|S_f\|=\|\iota\circ S_f\|$ or, equivalently, $f\in\widehat{\B}_\I(\D,X)$ with $\rho_\B(f)=\rho_\B(\iota\circ f)$.

(S) Let $f\in\widehat{\B}(\D,X)$ and let $\pi\colon\D\to\D$ be a holomorphic function with $\pi(0)=0$ such that $\widehat{\pi}\in\L(\G(\D),\G(\D))$ is a metric surjection. Assume that $f\circ\pi\in\widehat{\B}_\I(\D,X)$. Since $S_f\circ\widehat{\pi}=S_{f\circ\pi}\in\I(\G(\D),X)$ and the operator ideal $[\I,\|\cdot\|]$ is surjective by \cite[4.7.13]{Pie-80}, we have that $S_f\in\I(\G(\D),X)$ with $\|S_f\|=\|S_f\circ\widehat{\pi}\|$, hence $f\in\widehat{\B}_\I(\D,X)$ with $\rho_\B(f)=\rho_\B(f\circ\pi)$.
\end{proof}

Using Theorem \ref{prop2.1}, a proof similar to that of Proposition \ref{ideal} yields the following.

\begin{proposition}\label{ideal-00}
$[\widehat{\B}_\F,\rho_\B]$ is a surjective and injective normed normalized Bloch ideal.$\hfill\Box$
\end{proposition}


\subsection{Transposition}\label{section 5}

We state the analogues for Bloch mappings of the results due to Schauder, Gantmacher and Nakamura on the compactness and weak compactness of the adjoint of a bounded linear operator between Banach spaces. 

The first step is to present a Bloch version of the concept of adjoint operator. Let $f\in\widehat{\B}(\D,X)$. Given $x^*\in X^*$, we have $(x^*\circ f)(0)=0$ and 
\begin{align*}
(1-|z|^2)\|(x^*\circ f)'(z)\|&=(1-|z|^2)\|x^*(f'(z))\|\\
                                          &\leq (1-|z|^2)\|x^*\|\|f'(z)\|\leq \rho_\B(f)\|x^*\|
\end{align*}
for all $z\in\D$. Hence $x^*\circ f\in\widehat{\B}(\D)$ with $\rho_\B(x^*\circ f)\leq \rho_\B(f)\|x^*\|$. This justifies the following.

\begin{definition}
Given $f\in\widehat{\B}(\D,X)$, the \textit{Bloch transpose} of $f$ is the mapping $f^t\colon X^*\to\widehat{\B}(\D)$, given by 
$$
f^t(x^*)=x^*\circ f\qquad (x^*\in X^*).
$$
\end{definition}

Clearly, $f^t$ is linear and continuous with $\|f^t\|\leq \rho_\B(f)$. In fact, $\|f^t\|=\rho_\B(f)$. Indeed, for $0<\varepsilon<\rho_\B(f)$, take $z\in\D$ such that $(1-|z|^2)\|f'(z)\|>\rho_\B(f)-\varepsilon$. By the Hahn--Banach Theorem, there exists $y^*\in X^*$ with $\|y^*\|=1$ such that $\left|y^*(f'(z))\right|=\|f'(z)\|$. We have 
\begin{align*}
\|f^t\|&\geq\sup_{x^*\neq 0}\frac{\rho_\B(f^t(x^*))}{\|x^*\|}
\geq\frac{\rho_\B(y^*\circ f)}{\|y^*\|}\\
&\geq (1-|z|^2)\left|y^*(f'(z))\right|=(1-|z|^2)\|f'(z)\|>\rho_\B(f)-\varepsilon .
\end{align*}
Letting $\varepsilon\to 0$, one obtains $\|f^t\|\geq \rho_\B(f)$, as desired. Finally, note that  
\begin{align*}
(\Lambda\circ f^t)(x^*)(\gamma_z)&=\Lambda(f^t(x^*))(\gamma_z)=\Lambda(x^*\circ f)(\gamma_z)\\
                                 &=(x^*\circ f)'(z)=x^*(f'(z))\\
                                 &=x^*(S_f(\gamma_z))=(S_f)^*(x^*)(\gamma_z)
\end{align*}
for all $x^*\in X^*$ and $z\in\D$, where $(S_f)^*\colon X^*\to\G(\D)^*$ is the adjoint operator of $S_f$. Since $\G(\D)=\overline{\lin}(\Gamma(\D))$, we deduce that $\Lambda\circ f^t=(S_f)^*$ and thus $f^t=\Lambda^{-1}\circ(S_f)^*$. So we have:

\begin{proposition}\label{prop-A}
If $f\in\widehat{\B}(\D,X)$, then $f^t\in\L(X^*,\widehat{\B}(\D))$ with $\|f^t\|=\rho_\B(f)$ and $f^t=\Lambda^{-1}\circ(S_f)^*$. $\hfill\Box$
\end{proposition}

We next see that the mapping $f\mapsto f^t$ identifies $\widehat{\B}(\D,X)$ with the subspace of $\L(X^*,\widehat{\B}(\D))$ formed by all weak*-to-weak* continuous linear operators from $X^*$ into $\widehat{\B}(\D)$ (see \cite[Corollary 3.1.12]{Meg-98}).

\begin{proposition}\label{teo-4-1}
$f\mapsto f^t$ is an isometric isomorphism from $\widehat{\B}(\D,X)$ onto $\L((X^*,w^*);(\widehat{\B}(\D),w^*))$.
\end{proposition}

\begin{proof}
Let $f\in\widehat{\B}(\D,X)$. Hence $f^t=\Lambda^{-1}\circ (S_f)^*\in\L((X^*,w^*);(\widehat{\B}(\D),w^*))$ by Theorem \ref{teo-3} and \cite[Theorem 3.1.11]{Meg-98}.  
We have $\|f^t\|=\rho_\B(f)$ by Proposition \ref{prop-A}. It remains to show the surjectivity of the mapping in the statement. Take $T\in\L((X^*,w^*);(\widehat{\B}(\D),w^*))$. Then the mapping $\Lambda\circ T$ belongs to $\L((X^*,w^*);(\G(\D)^*,w^*))$ 
and, by \cite[Theorem 3.1.11]{Meg-98}, there is a $S\in\L(\G(\D),X)$ such that $S^*=\Lambda\circ T$. By Corollary \ref{cor-3-new}, there exists $f\in\widehat{\B}(\D,X)$ such that $S_f=S$. Hence $T=\Lambda^{-1}\circ(S_f)^*=f^t$, as desired. 
\end{proof}

The equivalence $(1)\Leftrightarrow (2)$ in the next result is a version of Schauder Theorem for compact Bloch mappings. 

\begin{theorem}\label{teo-4-2}
Let $f\in\widehat{\B}(\D,X)$. The following assertions are equivalent:
\begin{enumerate}
  \item $f\colon\D\to X$ is compact Bloch.
	\item $f^t\colon X^*\to\widehat{\B}(\D)$ is compact.
	\item $f^t\colon X^*\to\widehat{\B}(\D)$ is bounded-weak*-to-norm continuous. 
	\item $f^t\colon X^*\to\widehat{\B}(\D)$ is compact and bounded-weak*-to-weak continuous. 
	\item $f^t\colon X^*\to\widehat{\B}(\D)$ is compact and weak*-to-weak continuous.
\end{enumerate}
\end{theorem}

\begin{proof}
$(1)\Leftrightarrow (2)$: applying Theorem \ref{teo-3-1}, the Schauder Theorem \cite[Theorem 3.4.15]{Meg-98} and \cite[Proposition 3.4.10]{Meg-98}, we have 
\begin{align*}
f\in\widehat{\B}_\K(\D,X)
&\Leftrightarrow S_f\in\K(\G(\D),X)\\
&\Leftrightarrow (S_f)^*\in\K(X^*,\G(\D)^*)\\
&\Leftrightarrow f^t=\Lambda^{-1}\circ(S_f)^*\in\K(X^*,\widehat{\B}(\D)).
\end{align*}
$(1)\Leftrightarrow (3)$: similarly, one obtains 
\begin{align*}
f\in\widehat{\B}_\K(\D,X)
&\Leftrightarrow S_f\in\K(\G(\D),X)\\
&\Leftrightarrow (S_f)^*\in\L((X^*,bw^*);\G(\D)^*)\\
&\Leftrightarrow f^t=\Lambda^{-1}\circ(S_f)^*\in\L((X^*,bw^*);\widehat{\B}(\D)),
\end{align*}
by Theorem \ref{teo-3-1} and \cite[Theorem 3.4.16]{Meg-98}. 

$(3)\Leftrightarrow (4)\Leftrightarrow (5)$: it follows directly from \cite[Proposition 3.1]{Kim-13}.
\end{proof}

We are now in position to identify $\widehat{\B}_\K(\D,X)$ with the subspace of $\L((X^*,w^*);(\G(\D),w^*))$ consisting of all bounded-weak*-to-norm continuous linear operators from $X^*$ into $\widehat{\B}(\D)$.

\begin{proposition}\label{cor-4-1}
$f\mapsto f^t$ is an isometric isomorphism from $\widehat{\B}_\K(\D,X)$ onto $\L((X^*,bw^*);\widehat{\B}(\D))$.
\end{proposition}

\begin{proof}
Let $f\in\widehat{\B}_\K(\D,X)$. Then $f^t\in\L((X^*,bw^*);\widehat{\B}(\D))$ by Theorem \ref{teo-4-2} and $\|f^t\|
=\rho_\B(f)$ by Proposition \ref{prop-A}. To prove the surjectivity, take $T\in\L((X^*,bw^*);\widehat{\B}(\D))$. Then $\Lambda\circ T\in\L((X^*,bw^*);\G(\D)^*)$. If $Q_{\G(\D)}$ denotes the natural injection from $\G(\D)$ into $\G(\D)^{**}$, then $Q_{\G(\D)}(\gamma)\circ\Lambda\circ T\in\L((X^*,bw^*);\C)$ for all $\gamma\in\G(\D)$ and, by \cite[Theorem 2.7.8]{Meg-98}, $Q_{\G(\D)}(\gamma)\circ\Lambda\circ T\in\L((X^*,w^*);\C)$ for all $\gamma\in\G(\D)$, that is, $\Lambda\circ T\in\L((X^*,w^*);(\G(\D)^*,w^*))$ by \cite[Corollary 2.4.5]{Meg-98}. Hence $\Lambda\circ T=S^*$ for some $S\in\L(\G(\D),X)$ by \cite[Theorem 3.1.11]{Meg-98}. Note that $S^*\in\L((X^*,bw^*);\G(\D)^*)$ and this means that $S\in\K(\G(\D),X)$ by \cite[Theorem 3.4.16]{Meg-98}. Now, $S=S_f$ for some $f\in\widehat{\B}_{\K}(\D,X)$ by Theorem \ref{teo-3-1}. Finally, we have $T=\Lambda^{-1}\circ S^*=\Lambda^{-1}\circ (S_f)^*=f^t$. 
\end{proof}

The following result provides both versions of Gantmacher and Gantmacher--Nakamura Theorems for weakly compact Bloch mappings. 

\begin{theorem}\label{teo-4-3}
Let $f\in\widehat{\B}(\D,X)$. The following are equivalent:
\begin{enumerate}
  \item $f\colon\D\to X$ is weakly compact Bloch.
	\item $f^t\colon X^*\to\widehat{\B}(\D)$ is weakly compact.
	\item $f^t\colon X^*\to\widehat{\B}(\D)$ is weak*-to-weak continuous.
\end{enumerate}
\end{theorem}

\begin{proof}
$(1)\Leftrightarrow (2)$: we have  
\begin{align*}
f\in\widehat{\B}_\W(\D,X)
&\Leftrightarrow S_f\in\W(\G(\D),X)\\
&\Leftrightarrow (S_f)^*\in\W(X^*,\G(\D)^*)\\
&\Leftrightarrow f^t=\Lambda^{-1}\circ(S_f)^*\in\W(X^*,\widehat{\B}(\D)),
\end{align*}
by Theorem \ref{teo-3-2}, the Gantmacher Theorem \cite[Theorem 3.5.13]{Meg-98} and \cite[Proposition 3.5.11]{Meg-98}.

$(1)\Leftrightarrow (3)$: similarly, 
\begin{align*}
f\in\widehat{\B}_\W(\D,X)
&\Leftrightarrow S_f\in\W(\G(\D),X)\\
&\Leftrightarrow (S_f)^*\in\L((X^*,w^*);(\G(\D)^*,w))\\
&\Leftrightarrow f^t=\Lambda^{-1}\circ(S_f)^*\in\L((X^*,w^*);(\widehat{\B}(\D),w))
\end{align*}
by Theorem \ref{teo-3-2}, the Gantmacher--Nakamura Theorem \cite[Theorem 3.5.14]{Meg-98} and \cite[Corollary 2.5.12]{Meg-98}.
\end{proof}

We next identify $\widehat{\B}_\W(\D,X)$ with the subspace of $\L((X^*,w^*);(\G(\D),w^*))$ formed by all weak*-to-weak continuous linear operators from $X^*$ into $\widehat{\B}(\D)$.

\begin{proposition}\label{cor-4-1b}
$f\mapsto f^t$ is an isometric isomorphism from $\widehat{\B}_{\W}(\D,X)$ onto $\L((X^*,w^*);(\widehat{\B}(\D),w))$.
\end{proposition}

\begin{proof}
In view of Theorem \ref{teo-4-3} and Proposition \ref{prop-A}, we only need to show that the mapping in the statement is surjective. Let $T\in\L((X^*,w^*);(\widehat{\B}(\D),w))$. Then $\Lambda\circ T\in\L((X^*,w^*);(\G(\D)^*,w))$ by \cite[Theorem 2.5.11]{Meg-98}, and this last set is contained in the space $\L((X^*,w^*);(\G(\D)^*,w^*))$. 
It follows that $\Lambda\circ T=S^*$ for some $S\in\L(\G(\D),X)$ by \cite[Theorem 3.1.11]{Meg-98}. Hence $S^*\in\L((X^*,w^*);(\G(\D)^*,w))$ and, by the Gantmacher--Nakamura Theorem, $S\in\W(\G(\D),X)$. Now, $S=S_f$ for some $f\in\widehat{\B}_{\W}(\D,X)$ by Corollary \ref{teo-3-2}. Finally, $T=\Lambda^{-1}\circ S^*=\Lambda^{-1}\circ (S_f)^*=f^t$, as desired.
\end{proof}

We conclude by identifying the little Bloch space $\widehat{\B}_0(\D,X)$ with the space of bounded-weak*-to-norm continuous linear operators from $X^*$ to $\widehat{\B}_0(\D)$. In its proof we will apply a criterion for compactness in $\widehat{\B}_0(\D)$ due to Madigan and Matheson \cite{MadMat-95}.

\begin{theorem}\label{teo-4-4}
$f\mapsto f^t$ is an isometric isomorphism from $\widehat{\B}_0(\D,X)$ onto $\L((X^*,bw^*);\widehat{\B}_0(\D))$.
\end{theorem}

\begin{proof}
Let $f\in\widehat{\B}_0(\D,X)$ and let us prove that $f^t\in\L((X^*,bw^*);\widehat{\B}_0(\D))$. We have
$$
(1-|z|^2)\|f'(z)\|=(1-|z|^2)\sup_{x^*\in B_{X^*}}\left|x^*(f'(z))\right|=\sup_{x^*\in B_{X^*}}(1-|z|^2)\left|(x^*\circ f)'(z)\right|
$$
for every $z\in\D$. For each $x^*\in B_{X^*}$, we deduce that $f^t(x^*)\in\widehat{\B}_0(\D)$ with $\rho_\B(f^t(x^*))\leq \rho_\B(f)$. Hence $f^t(B_{X^*})$ is a bounded subset of $\widehat{\B}_0(\D)$. Since we also have  
$$
\lim_{|z|\to 1}\sup_{x^*\in B_{X^*}}(1-|z|^2)\left|(x^*\circ f)'(z)\right|=\lim_{|z|\to 1}(1-|z|^2)\|f'(z)\|=0,
$$
Lemma 1 in \cite{MadMat-95} asserts that the set $f^t(B_{X^*})$ is relatively compact in $\widehat{\B}_0(\D)$, that is, $f^t$ is in $\K(X^*,\widehat{\B}_0(\D))$. Moreover, $f^t$ belongs to $\L((X^*,w^*);(\widehat{\B}(\D),w^*))$ by Proposition \ref{teo-4-1}. Hence $\Lambda\circ f^t\in\K(X^*,\G(\D)^*)\cap\L((X^*,w^*);(\G(\D)^*,w^*))$. Then there is an $S\in\L(\G(\D),X)$ for which $S^*=\Lambda\circ f^t$ by \cite[Theorem 3.1.11]{Meg-98}. Hence $S^*\in\K(X^*,\G(\D)^*)$ and, by the Schauder Theorem, $S\in\K(\G(\D),X)$. Hence $f^t=\Lambda^{-1}\circ S^*\in\L((X^*,bw^*);\widehat{\B}(\D))$ by \cite[Theorem 3.4.16]{Meg-98}. Since $f^t(X^*)\subseteq\widehat{\B}_0(\D)$, we conclude that $f^t\in\L((X^*,bw^*);\widehat{\B}_0(\D))$, as expected.

By Proposition \ref{prop-A}, the mapping of the statement is a linear isometry. It remains to show that it is onto. For this take $T\in\L((X^*,bw^*);\widehat{\B}_0(\D))$. By Proposition \ref{cor-4-1}, there exists $f\in\widehat{\B}_\K(\D,X)$ such that $T=f^t$. Let us prove that $f$ actually belongs to $\widehat{\B}_0(\D,X)$. Since $T\in\K(X^*,\widehat{\B}(\D))$ by Theorem \ref{teo-4-2} and $T(X^*)\subseteq\widehat{\B}_0(\D)$, we deduce that $T\in\K(X^*,\widehat{\B}_0(\D))$. Then, by applying again \cite[Lemma 1]{MadMat-95} we obtain 
\begin{align*}
\lim_{|z|\to 1}(1-|z|^2)\|f'(z)\|&=\lim_{|z|\to 1}\sup_{x^*\in B_{X^*}}(1-|z|^2)\left|x^*(f'(z))\right|\\
&=\lim_{|z|\to 1}\sup_{x^*\in B_{X^*}}(1-|z|^2)\left|(T(x^*))'(z)\right|=0,
\end{align*}
and so $f\in\widehat{\B}_0(\D,X)$.  
\end{proof}


\textbf{Acknowledgements.} We would like to sincerely thank the reviewer for detecting a serious inconsistency in an earlier version of this paper, and the handling editor for indicating how to resolve this inconsistency.

A talk about this paper was given by the first author at the meeting ``Workshop on Function Spaces and Related Topics. The 5th Niigata Seminar'', held in Department of Mathematics of Niigata University from April 3--8, 2023. He is deeply grateful to Keiko (especially) and Osamu Hatori, Takeshi Miura, Shiho Oi and Rumi Shindo Togashi for their infinite hospitality during his unforgettable stay in Japan. 

Research funded by Junta de Andaluc\'{\i}a grant FQM194. Moreover, the first author was partially supported by Ministerio de Ciencia e Innovaci\'{o}n grant PID2021-122126NB-C31 (funded by MICIU/AEI/10.13039/501100011033 and by ``ERDF A way of making Europe''); and also by the Center for Development and Transfer of Mathematical Research to Industry CDTIME (University of Almer\'ia). The second author was supported by an FPU predoctoral fellowship of the Spanish Ministry of Universities (FPU23/03235).

\end{document}